\documentclass[twoside,letterpaper,11pt,leqno]{amsart}

\usepackage{fancyvrb}

\usepackage{url}
\usepackage[plainpages=false,pdfpagelabels,pdfencoding=auto,psdextra,colorlinks=true,urlcolor=black,linkcolor=black,citecolor=blue,hidelinks]{hyperref}

\usepackage[tmargin=1in, bmargin=1in, lmargin=1.49in, rmargin=1.49in]{geometry}

\usepackage{amsopn}
\usepackage{amsmath}
\usepackage{amsfonts}
\usepackage{amssymb}
\usepackage{graphicx}
\usepackage{float}
\usepackage{mathtools}
\usepackage[inline]{enumitem}
\usepackage{bm}
\usepackage[caption=false]{subfig}
\usepackage{multirow}
\usepackage[capitalize,nameinlink]{cleveref}
\usepackage{tikz}
\usepackage{tikz-qtree}
\usepackage{tikz-cd}
\usepackage{algorithm}
\usepackage{algorithmic}

\usepackage{microtype}

\let\Pi\varPi

\input{defines.texx}

\allowdisplaybreaks[4]

\newtheorem{theorem}{Theorem}

\theoremstyle{remark}
\newtheorem{remark}[theorem]{Remark}

\numberwithin{theorem}{section}
\numberwithin{equation}{section}

\numberwithin{algorithm}{section}

\numberwithin{figure}{section}
\numberwithin{table}{section}

\title[Parallel AMGe for $H(\ucurl)$ and $H(\div)$ Using ParELAG]{Parallel Element-based Algebraic Multigrid for $H(\ucurl)$ and $H(\div)$ Problems Using the ParELAG Library}
\author{Delyan Z. Kalchev} 
\address{Center for Applied Scientific Computing, Lawrence Livermore National Laboratory, P.O. Box 808, L-561, Livermore, CA 94551, USA.}
\email{kalchev1@llnl.gov}
\author{Panayot S. Vassilevski} 
\address{Department of Mathematics and Statistics, Portland State University, Portland, OR 97207, USA, and  Center for Applied Scientific Computing, Lawrence Livermore National Laboratory, P.O. Box 808, L-561, Livermore, CA 94551, USA.}
\email{panayot@pdx.edu, vassilevski1@llnl.gov}
\author{Umberto Villa} 
\address{Oden Institute for Computational Engineering \& Sciences, The University of Texas at Austin, Austin, TX 78712, USA.}
\email{uvilla@oden.utexas.edu}

\thanks{This work was performed under the auspices of the U.S. Department of Energy by Lawrence Livermore National Laboratory under contract DE-AC52-07NA27344 (LLNL-JRNL-824368).}
\thanks{The work of the second author was partially supported by NSF under grant DMS-1619640.}

\VerbatimFootnotes

\begin{document}

\begin{abstract}
This paper presents the use of element-based algebraic multigrid (AMGe) hierarchies, implemented in the ParELAG (Parallel Element Agglomeration Algebraic Multigrid Upscaling and Solvers) library, to produce multilevel preconditioners and solvers for $H(\ucurl)$ and $H(\div)$ formulations. ParELAG constructs hierarchies of compatible nested spaces, forming an exact de Rham sequence on each level. This allows the application of hybrid smoothers on all levels and AMS (Auxiliary-space Maxwell Solver) or ADS (Auxiliary-space Divergence Solver) on the coarsest levels, obtaining complete multigrid cycles. Numerical results are presented, showing the parallel performance of the proposed methods. As a part of the exposition, this paper demonstrates some of the capabilities of ParELAG and outlines some of the components and procedures within the library.


\end{abstract}

\maketitle

\section{Introduction}

Partial differential equation (PDE) models involving the $\ucurl$ (rotation) and $\div$ (divergence) operators often arise in the numerical simulation of physical phenomena and engineering systems. For example, this includes models of electromagnetism using Maxwell equations (possibly as a part of larger multiphysics codes) \cite{MonkFEM,2007LEforMHD,2003MagneticDiffusion}, mixed finite element methods for second-order elliptic equations \cite{BoffiMFE} and coupled systems \cite{2013Brinkman,2012VecLaplace}, first-order system least-squares (FOSLS) finite element methods \cite{1994FOSLS1,1996LSnonSA,2013ConstrFOSLS,AdlerVassilevski2014}, certain formulations of the Stokes and Navier-Stokes equations \cite{2010MixedStokes,2010MixedNS,1997NS}, and radiation transport simulations \cite{Brunner2002}.  

Among other discretization techniques, the finite element method (FEM) is a particularly appealing technique for problems defined on complex geometries due to its ability to handle unstructured meshes. One major challenge in solving linear systems arising from (FEM) discretizations of $H(\ucurl)$ and $H(\div)$ forms (i.e., symmetric problems involving the $\ucurl$ or $\div$ operators) is the large null spaces of the curl and divergence. A variety of approaches have been developed, including geometric and algebraic multigrid (AMG) \cite{1998Hiptmair,2002HcurlAMG,2003HcurlAMG,2006HcurlAMG,1997Hiptmair,2000HdivHcurlMultigrid,2004GenAMG,2008AMGkLaplacian,1992MixedMultigrid,1997HdivPecond,2009OptimalMG,VassilevskiMG}, static condensation and hybridization \cite{2019HybridizationHdiv}, and domain decomposition methods \cite{2018HdivBDDC,2016PETScBDDC,2016RobustnessBDDC,2002HcurlSchwartz,2000SchwarzHcurlHdiv}. Other techniques \cite{2006HcurlAux,2008HcurlAux} are based on reformulation of the governing equations that are attuned to (geometric) multigrid solvers and preconditioners.
A fundamental contribution to the development of multilevel methods for $H(\ucurl)$ and $H(\div)$ problems is the work of Hiptmair and Xu \cite{2007AuxiliarySpace}, which proposed auxiliary space preconditioners employing stable regular decompositions; see also \cite{2008AuxiliarySpace}. Based on those ideas quite successful parallel $H(\ucurl)$ and $H(\div)$ solvers were developed as a part of HYPRE \cite{hypre}: AMS \cite{2009AMS} and ADS \cite{2012ADS}. Furthermore, recent developments of generic auxiliary space preconditioners \cite{2020PreconditionerIP,2020PreconditionerMortar}, utilizing nonconforming reformulations and static condensation, can potentially be employed to implement efficient preconditioners for $H(\ucurl)$ and $H(\div)$ problems.

This paper describes and demonstrates the utilization of an element-based AMG (AMGe) approach (see, e.g., \cite{2001AMGe,vass_sparse_matrix_topology,VassilevskiMG}) for preconditioning conforming discrete $H(\ucurl)$ and $H(\div)$ formulations. While AMGe methods were originally developed in the context of symmetric positive definite (SPD) systems coming from $H^1$-conforming formulations \cite{2003AMGe,2008AMGe}, they demonstrate the capacity for a broader applicability. An important role in the construction of the AMGe multilevel methods is played by the \emph{de Rham complex}, a sequence of Sobolev spaces corresponding the domain and ranges of a chain of \emph{exterior derivatives}. For example, the sequence $H^1 \rightarrow H(\ucurl) \rightarrow H(\div) \rightarrow L^2$ form the three-dimensional de Rham complex, corresponding to the chain of $\grad$, $\ucurl$, $\div$ differential operators (exterior derivatives). For more detail regarding this elegant tool in the theory of finite elements, we refer to \cite{2002FEelectromagDeRham,2010DeRhamReview}. Namely, a discrete version of the sequence of conforming finite element spaces, maintaining the \emph{exactness} and \emph{commutativity} properties, delivers numerical stability (inf-sup compatibility) for a variety of mixed finite element methods. These ideas are used and investigated in \cite{2014CoarseDeRham,2012AMGeRT,2008deRhamAMGe} for the construction of multilevel element-based algebraic hierarchies of de Rham sequences of spaces. The last constitutes the foundation for the current work expounded in this paper.

A fundamental idea in AMGe, as presented in this work and implemented in ParELAG \cite{parelag}, is the element-based construction of coarse levels that structurally resemble (fine) geometric levels composed of standard finite elements. This involves the identification of coarse meshes with properly established \emph{coarse topologies} in the form of relations between coarse elements and coarse lower-dimensional mesh entities (facets, edges, and vertices), similarly to geometric levels. Consequently, utilizing the coarse topology, each coarse space is built via independent local coarse-element-by-coarse-element computations, whose combined effect is a conforming global coarse space. The independence of the local work makes the construction of AMGe hierarchies naturally attuned to parallel computing.

ParELAG is a parallel library that builds hierarchies of stable sequences of discrete spaces with approximation properties, to be utilized typically as discretization tools for numerical upscaling \cite{2011UpscalingAMG} of mixed finite element formulations. It also provides a set of respective preconditioners and solvers that can be used for solving the resulting problems or building composite solvers for more complex problems. ParELAG has been successfully applied, e.g., in upscaling for reservoir modeling \cite{2017ParELAGReservoir} and multilevel Monte Carlo simulations \cite{2018PDESampler,2017PDESampler,2021PosteriorMultilevel,2021HierarchicalMLMCMC}.

This paper discusses the construction of multilevel solvers for $H(\ucurl)$ and $H(\div)$ problems, using the hierarchies of spaces from ParELAG. To deliver a tidy and concrete presentation, ideas are conveyed for the three-dimensional case, while one can easily see how they would be applied in a two-dimensional setting. The availability of entire de Rham sequences, together with all necessary transfer operators, on all levels allows the utilization of \emph{hybrid (Hiptmair) smoothers} \cite{1998Hiptmair,VassilevskiMG} on all levels, as well as AMS and ADS on the coarsest levels, producing complete multigrid cycles. An outline of the overall methodology is presented and the parallel performance of the proposed solvers is shown in numerical examples. Finally, to further demonstrate ParELAG's capabilities and increase its visibility, a mini application within MFEM \cite{mfem}, a massively parallel widely used finite element library, has been developed.

The contributions and goals of this work (including vis-\`a-vis \cite{2008AMGe,2012AMGeRT,2014CoarseDeRham}) are: (1) to provide a unified presentation of AMGe methodologies for constructing coarse de Rham sequences by use of an exterior calculus formalism; (2) to describe the ParELAG library and exhibit its major capabilities; (3) to outline the generic multilevel AMGe methodology as implemented in the library, applicable to the entire de Rham sequence for coarse and fine finite element spaces of arbitrary order; and (4) to present a competitive, against state-of-the-art methods, novel $H(\ucurl)$ and $H(\div)$ AMGe solver with all its constitutive components (space hierarchies, smoothers, coarse solvers) implemented or invoked using the features of the core ParELAG library. From a theoretical perspective, a contribution of the paper is the use of the exterior calculus formalism to provide a unified dimension-independent presentation of the AMGe techniques of \cite{2008AMGe,2012AMGeRT,2014CoarseDeRham}. In particular, by use of the external derivative operator, the construction presented here for two and three-dimensional de Rham sequences can be extended to higher dimensions (e.g., to four-dimensional cases; see \cite{2018AuxiliaryDeRham}) or other complexes (such as the elasticity complex). From a software perspective, the utilization of the exterior calculus formalism allows for a drastic reduction of the volume of code in the implementation of the prolongator operators for both the two and thee-dimensional cases. Thus, our concise yet novel presentation of the methodologies serves to highlight unique features and capabilities of the ParELAG library. Finally, the numerical results serve the purpose of illustrating the potential and scalability of the technique in ParELAG, as well as the use of the library, which is the central topic of this paper and the MFEM mini application. Notably, the numerical results will also show that the AMGe solver with deep cycles demonstrate superior scalability, which motivates possible future work in fully exploiting this potential.

The outline of the remainder of the paper is as follows. The notation, the $H(\ucurl)$ and $H(\div)$ problems of interest, and an overview of de Rham sequences and their finite element discretization are presented in \Cref{sec:preliminaries}. \Cref{sec:HdeRham} is devoted to providing a succinct and unified (by use of exterior calculus formalism) presentation of the AMGe technique for the construction of a hierarchy of de Rham sequences in two or three space dimensions on agglomerated meshes, including the construction of prolongation operators, co-chain projectors, and coarse exterior derivatives. Those operators are useful for implementing the hybrid smoothers on all levels and the AMS and ADS coarse solvers, as described in \Cref{sec:smoothsolve}. Numerical results, demonstrating the parallel performance of the proposed methods, are in \Cref{sec:numerical}. The conclusions and a discussion of possible future directions are left for the last \Cref{sec:conclusions}.

\section{Preliminaries}
\label{sec:preliminaries}

This section presents the notation and function spaces used in this paper, as well as the formulation of the $H(\ucurl)$ and $H(\div)$ problems. It also provides an overview of finite element exterior calculus, including de Rham sequences of continuous and discrete (finite element) spaces. 

\subsection{\texorpdfstring{$H(\ucurl)$}{H(curl)} and \texorpdfstring{$H(\div)$}{H(div)} problems}

Let $\Omega \subset \bbR^3$ be a bounded contractible\footnote{Intuitively, a domain $\Omega$ can be continuously contracted to a point (i.e., is \emph{contractible}) if it has no holes or tunnels. Formally, a contractible domain is homotopy equivalent to a ball and to a single point.} domain with a Lipschitz-continuous boundary. Let $L^2(\Omega)$ and $[L^2(\Omega)]^3$ be the spaces of square integrable scalar and, respectively, vector functions. For $v \in L^2(\Omega)$ and $\uv \in [L^2(\Omega)]^3$, denote with $\lV v \rV^2_0 = (v, v)_0$ and $\lV \uv \rV^2_0 = (\uv, \uv)_0$ the norm induced by the corresponding inner product.

Let $D \in \{ \ugrad, \ucurl, \div, 0 \}$ denote an exterior derivative operator. Then, the notation $H(D\where \Omega)$ is used to denote the usual function spaces $H^1(\Omega) = \set{v \in L^2(\Omega)\where \ugrad v \in [L^2(\Omega)]^3}$, $H(\ucurl\where\Omega) = \set{ \uv \in [L^2(\Omega)]^3\where \ucurl \uv \in [L^2(\Omega)]^3}$, $H(\div\where\Omega) = \set{ \uv \in [L^2(\Omega)]^3\where \div \uv \in L^2(\Omega)}$, and $L^2(\Omega)$. Clearly, $H(D\where \Omega)$ is a Hilbert space endowed with the following inner product and induced norm:
\begin{equation}\label{eq:Dnorm}
\begin{split}
(\cdot, \cdot)_D &= (\cdot, \cdot)_0 + (D \cdot, D \cdot)_0,\\
\lV \cdot \rV_D &= \sqrt{(\cdot, \cdot)_D}.
\end{split}
\end{equation}
For ease of notation, when the domain is omitted, it is understood that $H(D) = H(D\where \Omega)$. 

Next, for $D\in\{\ucurl, \div\}$ consider the symmetric bilinear form 
\begin{equation}\label{eq:bfs}
\begin{alignedat}{4}
a_{D}(\uu, \uv) &= (\alpha\, D \uu, D \uv)_0 &&{}+{}&&(\beta\, \uu, \uv)_0\quad &&\text{ for }\; \uu, \uv \in H(D\where\Omega),
\end{alignedat}
\end{equation}
where $\alpha,\beta \in L^\infty(\Omega)$, $\alpha > 0$, $\beta > 0$. These bilinear forms, referred as $H(\ucurl)$ and $H(\div)$ forms hereafter, are positive definite and posses coefficient-dependent continuity in term of $\lV \cdot \rV_{D}$. If the coefficients are bounded away from zero, then the bilinear forms satisfy respective coefficient-dependent coercivity.

\begin{remark}
For simplicity, only the case $\beta > 0$ is studied here. However, the semi-definite case $\beta \ge 0$ is considered in \cite{2009AMS, 2012ADS}. In general, $\beta$ can be an essentially bounded symmetric positive (semi-)definite tensor. The bilinear forms are positive definite when $\beta$ is positive definite, generally semi-definite when $\beta$ is semi-definite, and coercivity depends on $\beta$ being uniformly (on $\Omega$) positive definite.
\end{remark}

By means of a Galerkin projection onto $H(D)$-conforming discrete finite element spaces, the discrete version of the bilinear form in \eqref{eq:bfs} can be represented by a symmetric positive (semi-)definite matrix. The goal of this paper is to construct multilevel preconditioners for linear systems with such matrices. The corresponding conforming finite element spaces, defined on a given fine mesh $\cT^h$, are denoted by $\cV^h(D) \subset H(D)$, for $D \in \{\ugrad, \ucurl, \div, 0\}$. They are spaces of, respectively, continuous piecewise polynomial Lagrangian (nodal), N{\'e}d{\'e}lec, Raviart--Thomas, and discontinuous piecewise polynomial finite elements \cite{BoffiMFE}. In the case of lowest order finite elements, the degrees of freedom (dofs) in the spaces are associated with mesh entities of increasing dimensionality, one dof per entity. Namely, these are, respectively, point values at vertices, tangential flow along edges, normal flux across facets, constant values in elements (sometimes referred to as cells).

\subsection{De Rham sequences of continuous and discrete spaces}
\label{ssec:deRhamSpaces}

\begin{figure}
\centering
\setlength\tabcolsep{1.5pt}
\begin{tabular}{ccccccc}
\begin{minipage}{0.19\textwidth}\includegraphics[width=\textwidth]{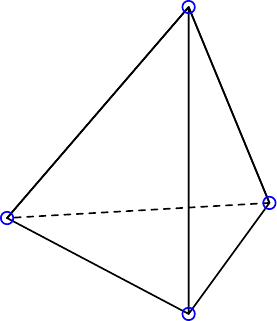}\end{minipage} &%
$\displaystyle\xrightarrow{\hphantom{.} D^h_1 \hphantom{.}}$ &%
\begin{minipage}{0.19\textwidth}\includegraphics[width=\textwidth]{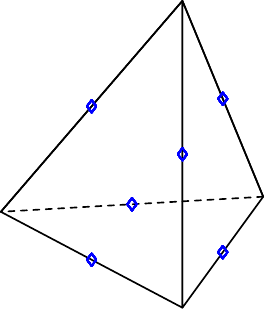}\end{minipage} &%
$\displaystyle\xrightarrow{\hphantom{.} D^h_2 \hphantom{.}}$ &%
\begin{minipage}{0.19\textwidth}\includegraphics[width=\textwidth]{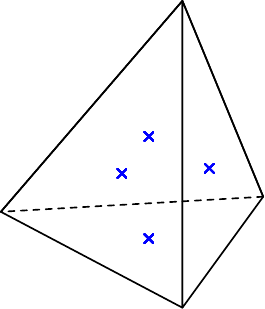}\end{minipage} &%
$\displaystyle\xrightarrow{\hphantom{.} D^h_3 \hphantom{.}}$ &%
\begin{minipage}{0.19\textwidth}\includegraphics[width=\textwidth]{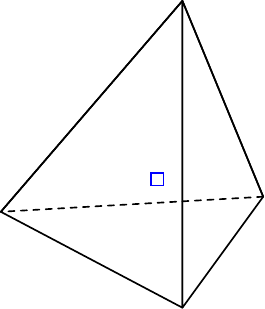}\end{minipage}
\end{tabular}
\caption{An illustration of the mapping between dofs in a tetrahedral element for the lowest order case.}\label{fig:dofslowestorder}
\end{figure}

Exterior calculus \cite{2010DeRhamReview} is a key tool in the stability and convergence analysis of finite element discretizations and solvers for the $H(D)$ form in \eqref{eq:bfs}. This is based on the de Rham complex of Sobolev spaces on $\Omega$ together with a respective \emph{subcomplex} of conforming finite element spaces
\begin{small}
\begin{equation}\label{eq:deRham}
\begin{tikzcd}[column sep=3em]
\bbR\arrow[r, "D_0 = \fI"] & H(D_1) \arrow[r, "D_1 = \ugrad"] \arrow[d, "\Pi^h_1"] & H(D_2) \arrow[r, "D_2 = \ucurl"] \arrow[d, "\Pi^h_2"] & H(D_3) \arrow[r, "D_3 = \div"] \arrow[d, "\Pi^h_3"] & H(D_4) \arrow[r, "D_4 = 0"] \arrow[d, "\Pi^h_4"] & \set{0}\\
\bbR\arrow[r, "D^h_0 = \fI"] & \cV^h(D_1) \arrow[r, "D^h_1"] & \cV^h(D_2) \arrow[r, "D^h_2"] & \cV^h(D_3) \arrow[r, "D^h_3"] & \cV^h(D_4) \arrow[r, "D^h_4 = 0"] & \set{0}.
\end{tikzcd}
\end{equation}
\end{small}%
Above, $\bbR$ represents the set of real numbers, $\fI$ is the injection operator mapping a real number to the corresponding constant function on $\Omega$, $\Pi^h_i: H(D_i) \to \cV^h(D_i)$ for $i=1,\dots,4$ are appropriate \emph{(cochain)} projection operators, 
$D_i: H(D_i) \to H(D_{i+1})$ for $i=1,\dots,3$ are differential operators (exterior derivatives) mapping between the Sobolev spaces, and $D^h_i: \cV^h(D_i) \to \cV^h(D_{i+1})$ are the corresponding discrete versions.

In the finite-dimensional setting, functions in $\cV^h(D_i)$ ($i=1,\dots,4$) can be identified with algebraic vectors collecting the coefficients in the respective finite element expansion (dofs). In what follows, $d^h_i=\dim(\cV^h(D_i))$ denotes the number of dofs of the space $\cV^h(D_i)$. Hence, $D^h_i$ for $i=1,\dots,3$ can be viewed as matrices in $\bbR^{d^h_{i+1}\times d^h_i}$ expressed in terms of the bases in $\cV^h(D_i)$. They can be assembled via an \emph{overwriting} finite element assembly\footnote{\emph{Overwriting} means that during the assembly the entries in the global matrix are overwritten by the values of the entries in the local matrices rather than accumulating (adding) them. See the implementation of \texttt{DiscreteOperator} in MFEM \cite{mfem}.
Since, as we will see in Section \ref{sec:HdeRham}, the coarse spaces generated by ParELAG are conforming across agglomerated entities interfaces, overwriting entries in the global matrix causes no alterations as values coming from neighboring entities are the same. The overwriting procedure is also in ParELAG to assemble the coarse exterior derivative operators $D^H_i$ and co-chain projectors $\Pi^H_i$, which are described below.
}
procedure from local, on elements (i.e., expressed in terms of shape functions), versions of the operators and their matrices. Note that, e.g., for the lowest order discretization (see \Cref{fig:dofslowestorder}), the operators $D^h_i$ ($i=1,\dots,3$) map from mesh entities of lower dimensionality to those of higher dimensionality, i.e., $\text{vertices}\to\text{edges}$, $\text{edges}\to\text{facets}$, and $\text{facets}\to\text{elements}$, respectively.

It is assumed that $\Pi^h_i$ for $i=1,\dots,4$ are bounded operators, i.e., $\lV \Pi^h_i \rV_{D_i} < \infty$, where $\lV \cdot \rV_{D_i}$ denotes the corresponding induced operator norm from \eqref{eq:Dnorm}. This holds for the considered finite element spaces and implies the quasi-optimality property $\lV u - \Pi^h_i u \rV_{D_i} \le \lV I - \Pi^h_i \rV_{D_i} \inf_{v^h \in \cV^h(D_i)} \lV u - v^h \rV_{D_i}$ for all $u \in H(D_i)$; see \cite{2010DeRhamReview}.

Furthermore, for $i=1,\dots,3$, let $\gamma_{\Gamma,i}$ denote the trace operator, that is the restriction of functions in $H(D_i \where \Omega)$ to $\Gamma \subset \partial \Omega$. Specifically, $\gamma_{\Gamma,1}$ is the usual trace operator mapping $H^1(\Omega) $ to $H^{1/2}(\Gamma)$; $\gamma_{\Gamma,2}$ restricts the tangential flow $\uv \times \un$ ($\uv \in H(D_2)$) to the surface $\Gamma$, and finally, $\gamma_{\Gamma,3}$ restricts the normal flux $\uv \cdot \un$ ($\uv \in H(D_3)$) to $\Gamma$. Here, $\un$ denotes the outward unit normal vector to $\partial \Omega$.

\begin{remark}
The diagram \eqref{eq:deRham} corresponds to the so-called \emph{natural} boundary conditions. It can be opportunely modified to the case of \emph{essential} boundary conditions (that is, the case of vanishing traces on $\partial \Omega$), as discussed in \cite{2010DeRhamReview}. 
\end{remark}

\subsubsection{Properties of the de Rham diagram}
Observe that the external derivative operators $D_i$ are such that $D_{i+1}D_i = 0$ (i.e., $\Range(D_i) \subset \Ker(D_{i+1})$) for $i=0,\dots,3$.

A de Rham sequence is called \emph{exact} if and only if 
\begin{equation}
\label{eq:exactness}
   \Range(D_i) = \Ker(D_{i+1}), 
\end{equation}
for $i=0,\dots,3$.
Exactness depends on the topological characteristics of $\Omega$. In particular, the connectivity of $\Omega$ is sufficient to demonstrate this property for $i=0\text{ and }3$, whereas it holds for $i=1$ using that $\Omega$ is simply-connected. The contractibility of $\Omega$ provides the property for $i=2$ as a consequence of Poincar{\'e}'s lemma; see, e.g., \cite{2002FEelectromagDeRham}.

The de Rham diagram \eqref{eq:deRham} is said to satisfy the \emph{commutativity} property if and only if
\begin{equation*}
D^h_i \circ \Pi^h_i = \Pi^h_{i+1} \circ D_i \quad\text{for } i=1,\dots,3.
\end{equation*}

Note that the commutativity property guarantees that, for a contractible domain $\Omega$, the exactness of the continuous de Rham sequence transfers to the discrete sequence \cite{2010DeRhamReview}.

The exactness of the continuous de Rham complex provides, e.g., stable decompositions (like the Helmholtz decomposition \cite{MonkFEM} and the so-called regular ones in \cite{2007AuxiliarySpace}), while the commutativity of \eqref{eq:deRham} and the exactness of the discrete subcomplex contribute to the inheritance of some important properties in the discrete setting, like the discrete stable decompositions in \cite{2007AuxiliarySpace} and the provision of the (inf-sup) stability of certain mixed finite element methods; see \cite{BoffiMFE,2002FEelectromagDeRham,2010DeRhamReview}. Such stability, together with the approximation properties of the discrete spaces, is a sufficient and necessary condition for the convergence of those mixed finite element methods. As is discussed in the following section, ParELAG builds de Rham sequences of coarse spaces satisfying the same exactness and commutativity properties of the fine sequence to ensure stability of coarse level discretizations, as well as approximation properties \cite{2014CoarseDeRham,2012AMGeRT,2008deRhamAMGe}.

\section{Overview of the multilevel de Rham sequence}
\label{sec:HdeRham}

\begin{figure}
\centering
\subfloat[][Agglomerates of fine elements]{\includegraphics[width=0.4\textwidth]{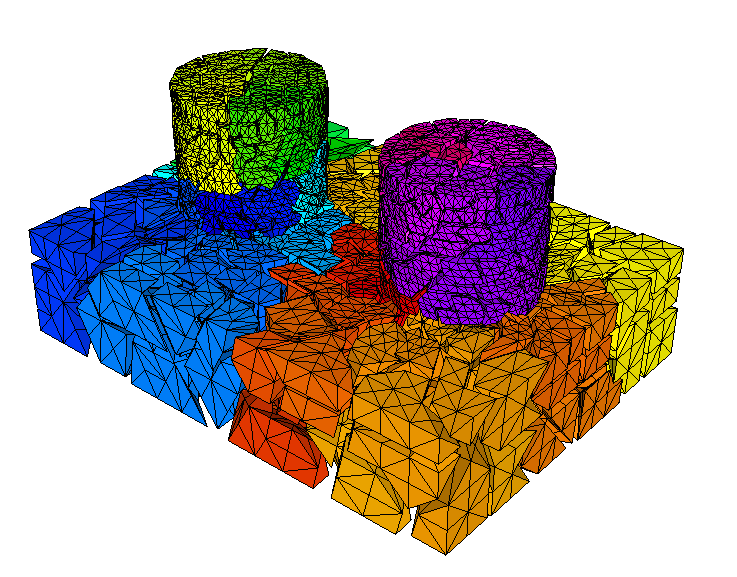}\label{fig:aggs}}\quad
\subfloat[][Agglomerates of agglomerates]{\includegraphics[width=0.4\textwidth]{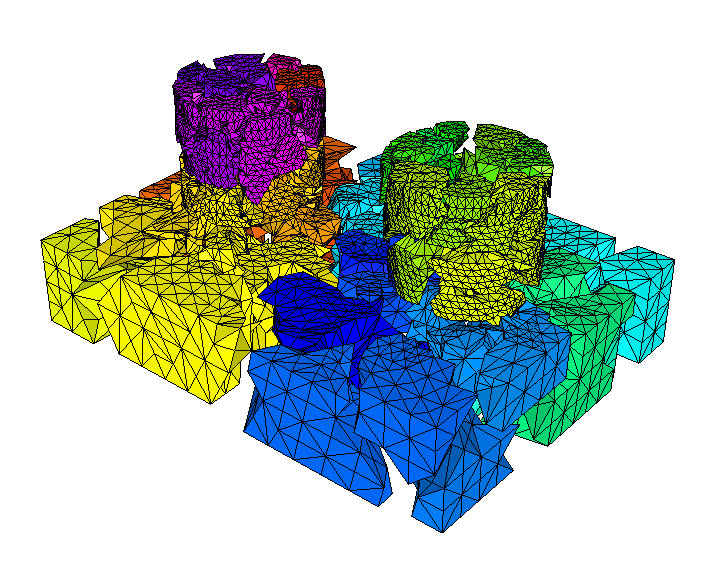}\label{fig:aggsaggs}}
\caption[]{Examples of agglomerates. (The gaps between agglomerates are for illustration.)}
\end{figure}

The fundamentals and notions associated with the coarsening of the de Rham sequence using AMGe techniques are now described. 
The key property of this approach, as articulated below, is that the de Rham sequence at each level of the hierarchy exhibits fine-like (geometric-like) finite element features. This means that coarse de 
Rham sequences are associated to generalized coarse meshes consisting of arbitrary shaped agglomerated entities (like elements, facets, edges, and vertices) and topological tables, and also maintain approximation properties as well as exactness (see \eqref{eq:Hexact}) and commutativity (see \eqref{eq:Hcomm}) properties, resulting in stable coarse de Rham complexes. This is a significant property of AMGe utilizing agglomeration of elements, which, together with the algebraic nature of the approach, allows the recursive application of the coarsening procedure. Doing so, ParELAG produces multilevel hierarchies of nested spaces forming exact and commutative de Rham complexes on all levels.

Starting with a given fine mesh the basic idea is to build a coarse mesh of coarse elements via \emph{agglomeration} and identify lower-dimensional coarse mesh entities, like coarse facets, edges, and vertices, together with their relationships, forming the coarse mesh topology. The mesh topology is needed for the construction of the de Rham sequences as it reflects a natural structure within the sequence of spaces and exterior derivatives. Note that the coarse de Rham sequence is formed in terms of the fine one as a sequence of subspaces, i.e., coarse basis functions are linear combinations of fine basis functions. The construction of coarse bases is purely algebraic and entails the obtainment of \emph{target traces}, or shortly \emph{targets}, and an \emph{extension} process. In particular, the extension process involves the solution of small local, on coarse entities, mixed finite element problems to produce the final coarse basis functions.

Specifically, consider the fine-level de Rham sequence (denoted with the superscript $h$) and the coarser level one (denoted with the superscript $H$), giving rise to the diagram
\begin{small}
\begin{equation}\label{eq:HdeRham}
\begin{tikzcd}[column sep=2.7em]
\bbR\arrow[r, "D^h_0 = \fI"] & \cV^h(D_1) \arrow[r, "D^h_1"] \arrow[d, "\Pi^H_1"] & \cV^h(D_2) \arrow[r, "D^h_2"] \arrow[d, "\Pi^H_2"] & \cV^h(D_3) \arrow[r, "D^h_3"] \arrow[d, "\Pi^H_3"] & \cV^h(D_4) \arrow[r, "D^h_4 = 0"] \arrow[d, "\Pi^H_4"] & \set{0}\\
\bbR\arrow[r, "D^H_0 = \fI"] & \cV^H(D_1) \arrow[r, "D^H_1"] & \cV^H(D_2) \arrow[r, "D^H_2"] & \cV^H(D_3) \arrow[r, "D^H_3"] & \cV^H(D_4) \arrow[r, "D^H_4 = 0"] & \set{0},
\end{tikzcd}
\end{equation}
\end{small}%
where $D_i^H: \cV^H(D_i)\to \cV^H(D_{i+1})$ ($i=1,\ldots,3$) and $\Pi^H_i: \cV^h(D_i) \to\cV^H(D_i)$ ($i=1,\ldots,4$) denote, respectively the coarse exterior derivative and co-chain projection matrices.

Assuming that the exactness \eqref{eq:exactness} of the continuous de Rham sequence holds, the coarsening procedure must guarantee the exactness property of the coarse sequence, that is,
\begin{equation}
\Range(D^H_i) = \Ker(D^H_{i+1})\quad \text{for } i = 1,\dots,3, \label{eq:Hexact}\\
\end{equation}
and the commutativity property
\begin{equation}
    D^H_i \circ \Pi^H_i = \Pi^H_{i+1} \circ D^h_i\quad \text{for } i = 1,\dots,3.\label{eq:Hcomm}
\end{equation}
These properties of the general approach to construct exact and commutative hierarchies via element-based algebraic multigrid techniques were proved in Pasciak and Vassilevski \cite{2008deRhamAMGe} (for the lowest order case) and in Lashuk and Vassilevski \cite{2014CoarseDeRham} (in the general case).

To this aim, ParELAG not only builds prolongation operators $P_i: \cV^H(D_i) \to\cV^h(D_i)$ (which allows to transfer information between levels), but also the coarse exterior derivative operators $D_i^H$ (which allows to define the hybrid smoothers in \Cref{ssec:hsmoothers}), and co-chain projectors $\Pi_i^H$ (which play a fundamental role in applying the auxiliary space AMG preconditioners to coarse problems as described in \Cref{ssec:AMSADS}). Moreover, ParELAG keeps track of element and facets attributes to ensure that material properties and essential boundary conditions can be properly applied to the discretized systems at every level of the hierarchy.

Note that the restriction, $P_i^T: \cV^h(D_i) \to\cV^H(D_i)$, and co-chain projection, $\Pi^H_i$, operators represent different actions and should not be confused; see \Cref{sssec:cprojcderiv}.

Finally, it is important to note that ParELAG need not necessarily build the entire sequence of spaces, if the application does not require it. Instead, ParELAG builds spaces of the de Rham sequence in reverse order, that is from $H(D_4) \equiv L^2$ to $H(D_i)$ for any $i \in [1,4]$. 

\subsection{Construction of a hierarchy of coarse meshes by agglomerating finer mesh entities}
\label{sssec:Hmesh}

The first step is the generation of a coarse mesh $\cT^H$, from the given fine one $\cT^h$, including all mesh entities: coarse elements, facets, edges, and vertices. The foundation of this is the construction of coarse elements as \emph{agglomerates} (or \emph{agglomerated elements}), which provide a non-overlapping partition of the fine elements; see \Cref{fig:aggs}. This is performed in a recursive manner to generate a hierarchy of nested meshes; see \Cref{fig:aggsaggs}. One customary way to achieve that is via partitioning (e.g., using METIS \cite{metis}) of the \emph{dual graph} of $\cT^h$---a graph whose nodes are the elements in $\cT^h$ and any two nodes are connected in the graph when the respective mesh elements share a facet. It is not difficult to generate the agglomerates as contiguous partitions in terms of the dual graph, e.g., using METIS or simply identifying the connected components of the partitioning after it is generated. Moreover, ParELAG provides additional tools that, via weighting the dual graph and further splitting of agglomerates, can help improve the topological properties of the coarse elements, which are relevant if $H(\ucurl)$ is utilized. 

\begin{figure}
\centering
\includegraphics[width=0.4\textwidth]{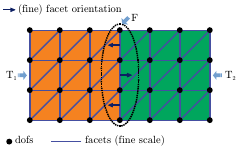}
\caption[]{A two-dimensional illustration of the designation of a coarse facet $F$ as a set of fine-scale facets, serving as an interface between agglomerates $T_1$ and $T_2$. Arrows illustrate the (global) orientation of the fine-scale facets, i.e., the orientation of the vectors normal to the facets.}\label{fig:facets}
\end{figure}

Using the partitioning of $\cT^h$ and viewing each agglomerate $T \in \cT^H$ as a collection of fine facets, an intersection procedure (see \cite[Section 1.9]{VassilevskiMG}) over these collections provides the coarse facets as sets of fine facets (see \Cref{fig:facets}), which can be consistently interpreted as interface surfaces between coarse elements. Further viewing the obtained coarse facets as collections of fine edges, their intersection identifies coarse edges as sets of fine edges. Finally, the intersection of coarse edges in terms of fine vertices identifies the coarse vertices.

Coarse facets and edges additionally carry information about the orientation of their constituting fine entities; see \Cref{fig:facets}. Such a set of fine-scale orientations for a coarse entity represents the orientation of that coarse entity. More precisely, these are $+1$ and $-1$ data entries in the agglomerated topology relating each coarse entity to its comprising fine-scale ones, respectively representing the preservation or the reversal of the original orientation of the fine entity within the coarse-scale one, so that each agglomerated entity has a consistent orientation. For example, a coarse facet $F$ has an associated vector of $+1$ and $-1$, denoted by $\bvarphi^F$, that based on the orientation of the constituting fine facets devises a consistent orientation for $F$, so that the normal vector to $F$ points everywhere from one of its adjacent agglomerates to the other, e.g., from $T_1$ to $T_2$ in \Cref{fig:facets}. Following the example in \Cref{fig:facets}, for the coarse facet $F$ to be oriented so that its normal vector consistently points from $T_1$ to $T_2$, one has $\bvarphi^F = [-1, +1, -1]$, indicating that the orientations of the first and third constituting fine facets are flipped while the orientation of the second one is preserved.

Finally, coarse entities associated with the domain boundary, or portions of it, are identified, thus allowing to track boundary dofs and apply boundary conditions on coarse levels. 

\subsubsection{ParELAG's implementation of agglomerated meshes hierarchy}

ParELAG contains a set of \emph{partitioner} classes, which generate an element partitioning on the current level that composes the agglomerated elements. For example, the class \verb|MFEMRefinedMeshPartitioner| constructs agglomerates in the form of geometric coarse elements by reverting previous refinements performed by MFEM, while the class \verb|MetisGraphPartitioner| invokes METIS internally. Furthermore, the coarse elements can be optionally made to conform to material (coefficient) interfaces by splitting agglomerates that cross such interfaces.

Instances of \verb|AgglomeratedTopology| collect coarse entities together with relationships between them in the form of a so called \emph{(agglomerated) topology} of $\cT^H$. Note that such a topology object in itself also represents a complex related to \eqref{eq:HdeRham}. It further contains relations between agglomerated entities and their constituting fine ones. Particularly, the \verb|AgglomeratedTopology| object on the finest level is obtained from the given mesh $\cT^h$ (i.e., using MFEM). Relationships between fine and coarse entities are stored in \verb|AgglomeratedEntity_Entity| tables within \verb|AgglomeratedTopology|. Each row in  \verb|AgglomeratedEntity_Entity| table corresponds to a coarse entity and the non-zeros entries in each row represent the finer grid entities that form the coarse entity. Furthermore, \verb|AgglomeratedEntity_Entity| tables for facets and edges also store orientation information.

Finally, having defined a topology on the current level, a new coarser agglomerated topology is generated by invoking the \verb|CoarsenLocalPartitioning()| member function of \verb|AgglomeratedTopology|, using the agglomerated elements produced on the current level by a partitioner class.

\subsection{The element-based construction of the components of the coarse sequences}
\label{sssec:ebasedP}

Here, we set the stage for the definition of the coarse basis functions in \Cref{ssec:Hbases}. We exhibit and motivate the fundamental utility of the process of building coarse bases by reviewing the element-based construction of the components that constitute the coarse sequences and coarse problems, once the coarse basis functions are obtained via the procedures in \Cref{ssec:Hbases}. This effectively reduces the considerations to the local agglomerate-by-agglomerate algebraic process of building locally-supported coarse bases presented in \Cref{ssec:Hbases}. Particularly, we review the coarse spaces formed via the choice of coarse bases, the construction of prolongation and coarse system matrices, including their local versions, as well as co-chain projectors and coarse exterior derivatives. In the end, we comment on the ParELAG classes implementing the multilevel de Rham sequences of spaces.

The procedures here are outlined to the minimal necessary extent following a concise constructive or algorithmic perspective and presented as implemented in ParELAG, since the goal is to describe the library and its capabilities. Therefore, the constructs are exposed and reviewed concisely for the purpose of presenting the software. Full details on the theoretical justification of the constructions explained in what follows are found in Pasciak and Vassilevski \cite{2008deRhamAMGe} (for the lowest order case) and in Lashuk and Vassilevski \cite{2014CoarseDeRham} (in the general case).

\subsubsection{Coarse spaces, prolongation matrices, and coarse system matrices}
\label{sssec:cbasescsystems}

\begin{figure}
\centering
\subfloat{\resizebox{0.8\textwidth}{!}{\input{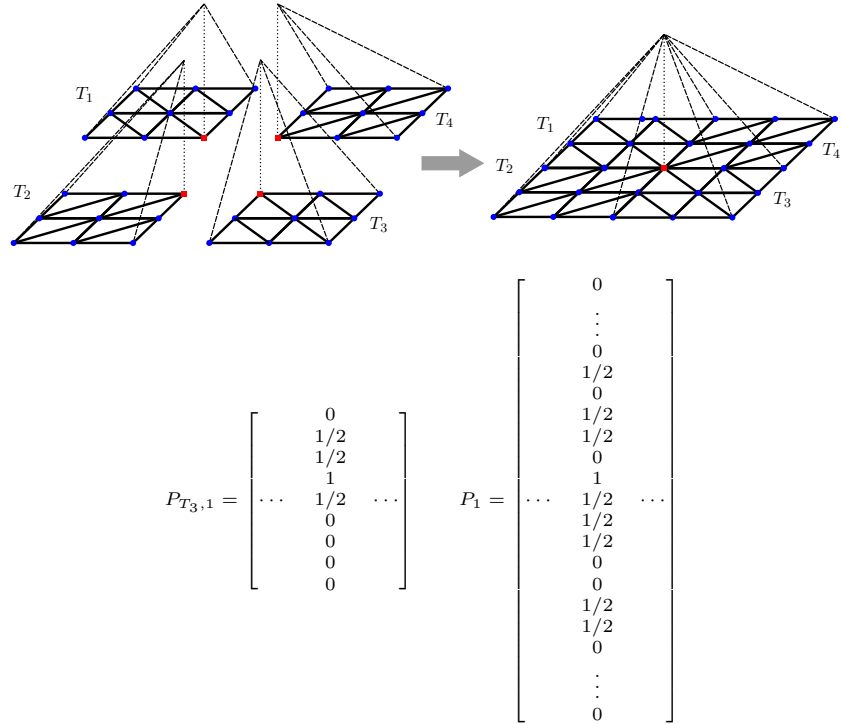}}}\\
\subfloat{\tiny
$
P_{T_3,1} = 
\begin{bmatrix}
&0&\\
&\halff&\\
&\halff&\\
&1&\\
\cdots&\halff&\cdots\\
&0&\\
&0&\\
&0&\\
&0&\\
\end{bmatrix}\quad\quad
P_{1} = 
\begin{bmatrix}
&0&\\
&\vdots&\\
&0&\\
&\halff&\\
&0&\\
&\halff&\\
&\halff&\\
&0&\\
&1&\\
\cdots&\halff&\cdots\\
&\halff&\\
&\halff&\\
&0&\\
&0&\\
&\halff&\\
&\halff&\\
&0&\\
&\vdots&\\
&0&\\
\end{bmatrix}
$
}
\caption[]{An illustration of coarse shape and basis functions and interpolation matrices. The example presents a piecewise linear $H^1$-conforming space on a two-dimensional mesh, where a single coarse basis function is supported on four agglomerates. The top left illustrates the separate shape functions that form the basis function. The final basis function (top right) is the result of conformingly combining the shape functions and the respective vectors defined on the fine dofs. The corresponding coarse dof is indicated by a small rectangle, while the fine dofs include the ones indicated by small circles and the rectangular one. The bottom shows a possible algebraic representation as column vectors of respective local (e.g., on $T_3$) and global interpolation matrices under some dof indexing.}\label{fig:shapebasis}
\end{figure}

The coarse spaces $\cV^H(D_i)$ for $i=1,\dots,4$ are obtained via the construction of coarse bases in \Cref{ssec:Hbases} as sets of algebraic vectors in terms of the respective dofs in $\cV^h(D_i)$, i.e., coarse basis functions are linear combinations of fine basis functions. These algebraic vectors constitute the columns of corresponding \emph{prolongation} matrices $P_i \in \bbR^{d^h_i \times d^H_i}$, $P_i: \cV^H(D_i) \to \cV^h(D_i)$, with full column ranks, where $d^H_i = \dim(\cV^H(D_i))$. Similarly to the fine level, the coarse basis functions are supported locally and built locally agglomerated entity by agglomerated entity (more details in \Cref{ssec:Hbases}). The $\cV^H(D_i)$-dofs are identified with the columns of $P_i$, i.e., with the respective coarse basis functions. Moreover, the $\cV^H(D_i)$-dofs associated with an agglomerate $T \in \cT^H$ are the ones whose basis functions have supports intersecting $T$, and the restrictions of those basis functions (i.e., the coarse \emph{shape functions}) on $T$ are precisely the restrictions of the respective algebraic vectors on the $\cV^h(D_i)$-dofs of $T$, which are the $\cV^h(D_i)$-dofs associated with the fine elements $\tau \in \cT^h$ such that $\tau \subset T$. Therefore, local-on-$T$ prolongation matrices $P_{T,i}$ can be defined and they are submatrices of $P_i$ on the respective dofs in $\cV^h(D_i)$ and $\cV^H(D_i)$ on $T$; see \Cref{fig:shapebasis}.

Next, consider the construction of system matrices at coarse levels of the AMGe hierarchy. Let $a_{ij}(\cdot,\cdot)$ be a bilinear form defined on $H(D_i)\times H(D_j)$ for some $i,j \in \{1,\dots,4\}$. By means of a Galerkin projection onto the finite element bases of the conforming discrete subspaces $\cV^h(D_i)$ and $\cV^h(D_j)$, the bilinear form is represented by a (global) matrix $A^h_{ij} \in \bbR^{d^h_j \times d^h_i}$ on the dofs in $\cV^h(D_i)$ and $\cV^h(D_j)$. That is, for every entry of $A^h_{ij}$ indexed $(l,k)$, it holds
\[
(A^h_{ij})_{lk} = a_{ij}(\phi^h_{i,k},\, \phi^h_{j,l})\quad\text{for }\;l = 1,\dots,d^h_j,\;k = 1,\dots,d^h_i,
\]
where $\set{\phi^h_{i,k}}_{k=1}^{d^h_i}$ denotes the basis of $\cV^h(D_i)$. This global matrix is obtained via a standard assembly from local element matrices $A^h_{\tau,ij}$ for the elements $\tau \in \cT^h$ formulated on the $\cV^h(D_i)$ and $\cV^h(D_j)$-dofs associated with $\tau$. The coarse matrices are produced by standard ``RAP'' procedures. Indeed, the representations of $a_{ij}(\cdot,\cdot)$ in terms of the bases of $\cV^H(D_i)$ and $\cV^H(D_j)$ is the matrix $A^H_{ij} = P_j^T A^h_{ij} P_i \in \bbR^{d^H_j \times d^H_i}$. Also, for $T\in\cT^H$, using a standard assembly locally with $A^h_{\tau,ij}$ for $\tau \subset T$, the local-on-$T$ fine-scale matrix $A^h_{T,ij}$ is obtained on the $\cV^h(D_i)$ and $\cV^h(D_j)$-dofs associated with $T$. Thus, $A^H_{T,ij} = (P_{T,j})^T A^h_{T,ij} P_{T,i}$ forms the coarse element matrices, which can produce $A^H_{ij}$ via a standard assembly. 

\subsubsection{Co-chain projectors and coarse exterior derivatives}
\label{sssec:cprojcderiv}

Here, we comment on the co-chain projection operators as well as on the coarse derivative operators. Co-chain projectors $\Pi_i^H$ are (right) inverses of the prolongation operators $P_i$, i.e., they satisfy
\begin{equation}
\label{eq:co-chain-proj-right-inverse}
 \Pi_i^H \circ P_i = I^H_i \text{ for } i = 1,\ldots,4,
\end{equation}
where $I^H_i$ denotes the identity operator in $\cV^H(D_i)$.
Their construction parallels that of the coarse basis functions in \Cref{ssec:Hbases}, starting from the projection of coarse dofs associated with local lower-dimensional entities and  moving towards higher-dimensional local entities. Specifically, projection operators are obtained via independent local agglomerate-by-agglomerate procedures. The local projection operators can be pieced together to form the global $\Pi^H_i$, since, by construction, these local projectors agree on dofs shared between coarse entities. The independent production of each such local projection operator involves the inversion of a small coarse-scale local mass matrix on the respective coarse entity. More details can be found in \cite{2014CoarseDeRham}; see also \cite{2008deRhamAMGe,2012AMGeRT}.

\begin{remark}
Note that the co-chain projector $\Pi_i^H$ is a left inverse of the prologation operator $P_i$ that satisfies two important properties: (1) it can be constructed locally (agglomerated entity by agglomerated entity) and (2) satisfies the commutative property \eqref{eq:Hcomm}. Thus, $\Pi_i^H$ differs from the restriction $P_i^T$, which may not satisfy \eqref{eq:co-chain-proj-right-inverse} or \eqref{eq:Hcomm}.
As such, $\Pi_i^H$ is a projection operator from fine to coarse finite element function spaces, while $P_i^T$ can be viewed as acting on functionals (i.e., mapping dual spaces). That is, $P_i^T$ provides restrictions suitable for Galerkin projections of linear or bilinear forms as elaborated in \Cref{sssec:cbasescsystems}. Specifically, $P_i^T: [\cV^h(D_i)]' \to [\cV^H(D_i)]'$, where $\ell^H = P_i^T \ell^h \in [\cV^H(D_i)]'$, for $\ell^h \in [\cV^h(D_i)]'$, acts like $\ell^H(v^H) = \ell^h(P_i v^H)$, for all $v^H \in \cV^H(D_i)$. That is, as the name suggests, $P_i^T$ restricts the action of $\ell^h$ from $\cV^h(D_i)$ to $\Range(P_i) = \cV^H(D_i)$.
\end{remark}

It is worth highlighting that co-chain projection operators do not merely constitute a theoretical tool, but they are needed for the construction of the coarse level auxiliary space preconditioners described in \Cref{ssec:AMSADS}. Moreover, they are needed in the implementation of multilevel Monte Carlo methodologies \cite{2018PDESampler,2017PDESampler} and efficient multilevel nonlinear solvers like FAS (full approximation scheme) \cite{2018AMGeFAS}. 

Next, consider the coarse exterior derivatives. The commutativity property \eqref{eq:Hcomm} implies that the coarse derivatives $D_i^H$ can be formally defined in terms of the prolongation $P_i$, cochain projector $\Pi_{i+1}^H$, and finer level exterior derivative $D_i^h$ as
\begin{equation*}
    D_i^H = \Pi_{i+1}^H \circ D_i^h \circ P_i,
\end{equation*}
as can be easily shown by (right-)multiplying \eqref{eq:Hcomm} by $P_i$ and using \eqref{eq:co-chain-proj-right-inverse}.

In ParELAG, the coarse exterior derivatives $D_i^H$ are constructed via independent local agglomerate-by-agglomerate procedures starting from lower-dimensional entities (agglomerated faces for $D_3^H$, agglomerated edges for $D_2^H$, and agglomerated vertices for $D_1^H$) to agglomerated elements $T$. The global exterior derivative operator $D^H_i$ is then constructed by an overwriting assembly of local-on-$T$ coarse exterior derivatives $D^H_{T,i}$ in a similar manner to the assembly of the finite element level operators discussed in \Cref{ssec:deRhamSpaces}. Note that, differently from matrices arising from discretizations of variational forms, at each level of the hierarchy the restrictions of the exterior derivative to a single entity or to the union of connected entities (agglomerate) are submatrices of $D^H_i$. 

\subsubsection{ParELAG's classes for AMGe hierarchies of de Rham sequences}

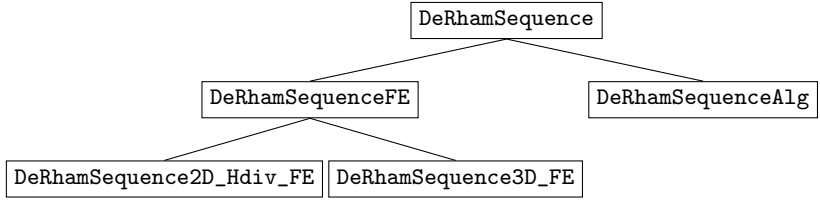
\begin{figure}
\centering
\footnotesize
\SaveVerb{DeRhamSequence2DHdivFE}|DeRhamSequence2D_Hdiv_FE|
\SaveVerb{DeRhamSequence3DFE}|DeRhamSequence3D_FE|
\begin{tikzpicture}[every tree node/.style=draw]
\tikzset{level distance=30pt}
\Tree [.{\Verb|DeRhamSequence|}
          [.{\Verb|DeRhamSequenceFE|}
              {\UseVerb{DeRhamSequence2DHdivFE}}
              {\UseVerb{DeRhamSequence3DFE}}
          ]
          [.{\Verb|DeRhamSequenceAlg|} ]
      ]
\end{tikzpicture}
\caption[]{An illustration of the hierarchy of ParELAG classes for de Rham sequences.}\label{fig:classes}
\end{figure}

ParELAG implements the construction and coarsening of de Rham sequences via the small hierarchy of classes depicted in \Cref{fig:classes}. The base class is \verb|DeRhamSequence|, which contains the main toolset necessary for constructing and working with de Rham sequences, including the procedures for building coarse spaces outlined in this paper. Two subclasses provide specialized methods for de Rham sequences on the finest (geometric) level (\verb|DeRhamSequenceFE|), which is produced employing MFEM, and \emph{algebraic} levels (\verb|DeRhamSequenceAlg|), which are coarse levels produced by ParELAG that are not associated with a given mesh (i.e., that are not \emph{geometric}). The class \verb|DeRhamSequenceFE| is further specialized to address special cases like the dimensionality of the domain.

\subsection{Coarse bases and the extension procedure}
\label{ssec:Hbases}

The abstract construction of coarse basis functions, which is applicable on all levels, is outlined now. The target traces and the extension process are discussed. A detailed presentation of a closely related procedure for coarse space construction can be found in \cite{2014CoarseDeRham}. After a short introduction, the procedure is described, followed by a schematic summary in \Cref{sssec:summary}.

The coarse mesh $\cT^H$ with all its entities (elements or agglomerates, facets, edges, and vertices) and its topology is available. Note that coarse basis functions are obtained in terms of the fine ones using an extension procedure involving the solution of local finite element problems, which translates into inverting local matrices. This uses information on the association between fine dofs and coarse entities, which is easily derived from the association between fine dofs and the fine entities that constitute each coarse entity. In the terminology of ParELAG, this is called \emph{dof agglomeration} (or \emph{aggregation}), implemented in the \verb|DofAgglomeration| class.

\begin{figure}
\centering
\captionsetup[subfloat]{font=tiny,width=0.3\textwidth}
\subfloat[][Coarse elements]{\includegraphics[width=0.3\textwidth]{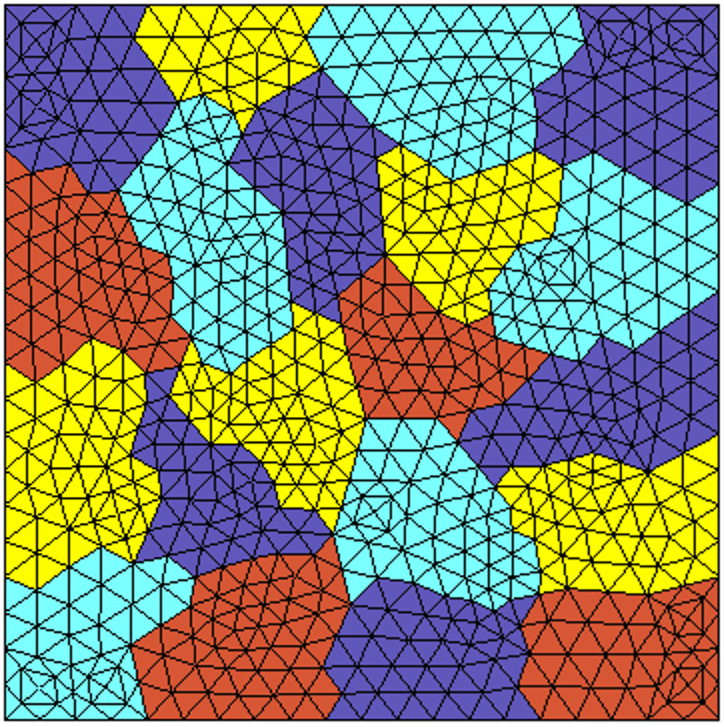}}\;
\subfloat[][$H^1$ boundary basis function]{\includegraphics[width=0.3\textwidth]{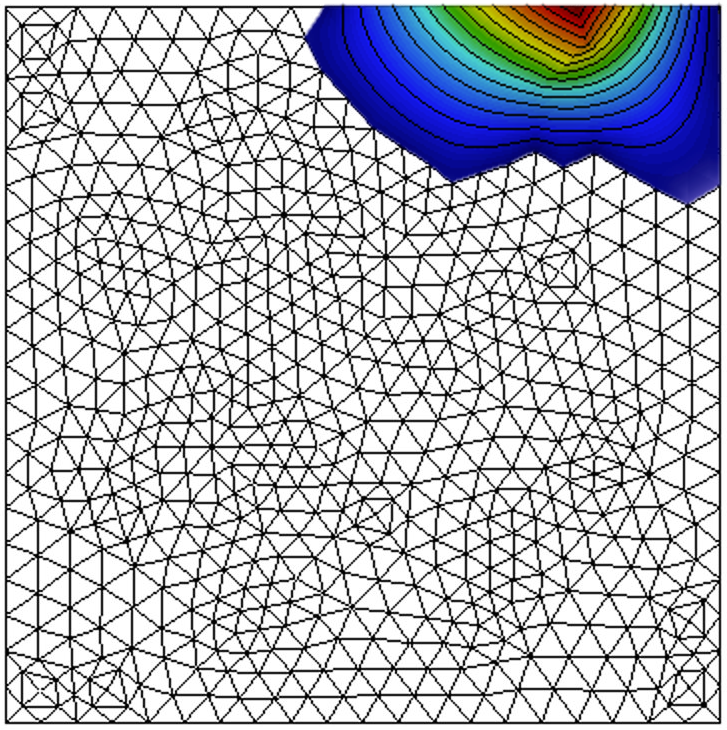}}\;
\subfloat[][$H^1$ interior basis function]{\includegraphics[width=0.3\textwidth]{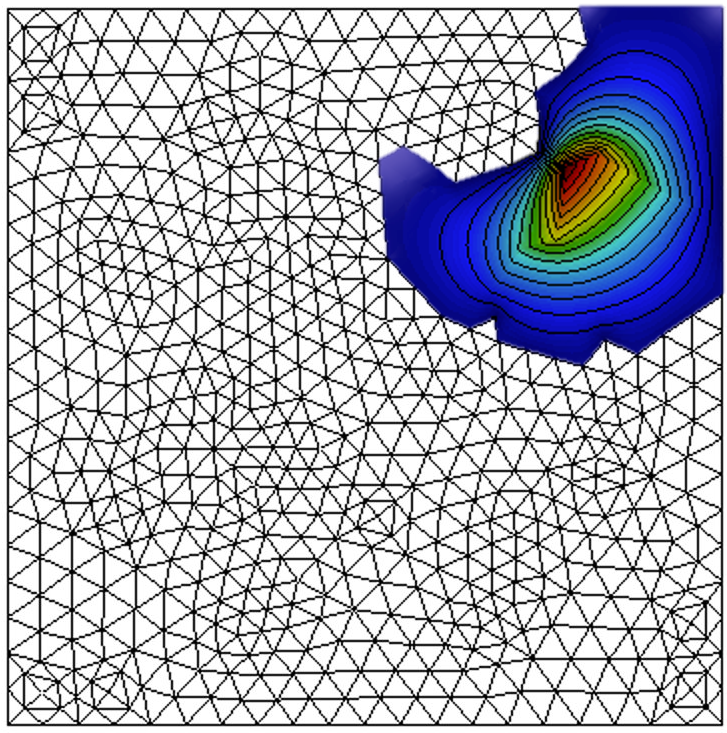}}\\
\subfloat[][$H^1$ bubble basis function]{\includegraphics[width=0.3\textwidth]{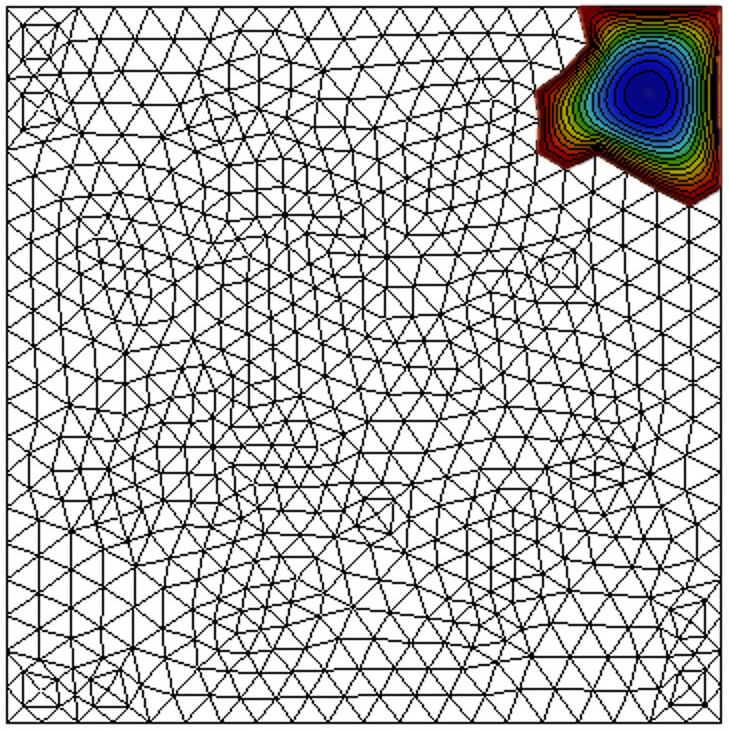}}\;
\subfloat[][$H(\div)$ boundary basis function]{\includegraphics[width=0.3\textwidth]{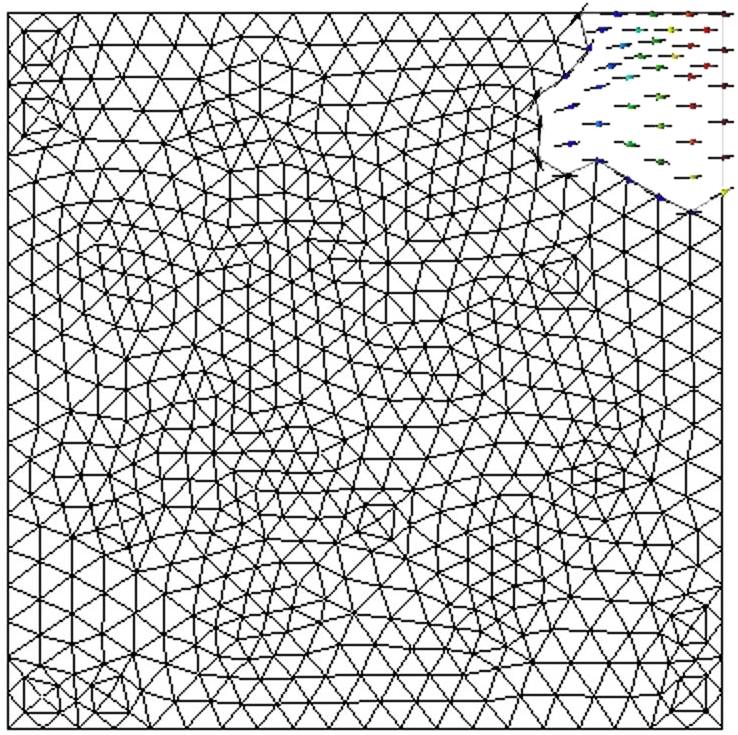}}\;
\subfloat[][$H(\div)$ interior basis function]{\includegraphics[width=0.3\textwidth]{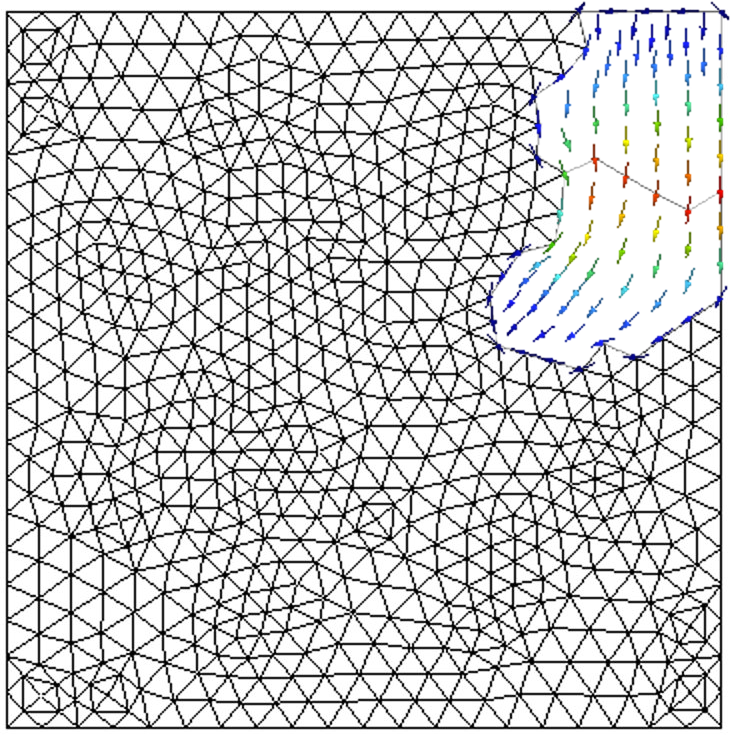}}\\
\subfloat[][$H(\div)$ bubble basis function]{\includegraphics[width=0.3\textwidth]{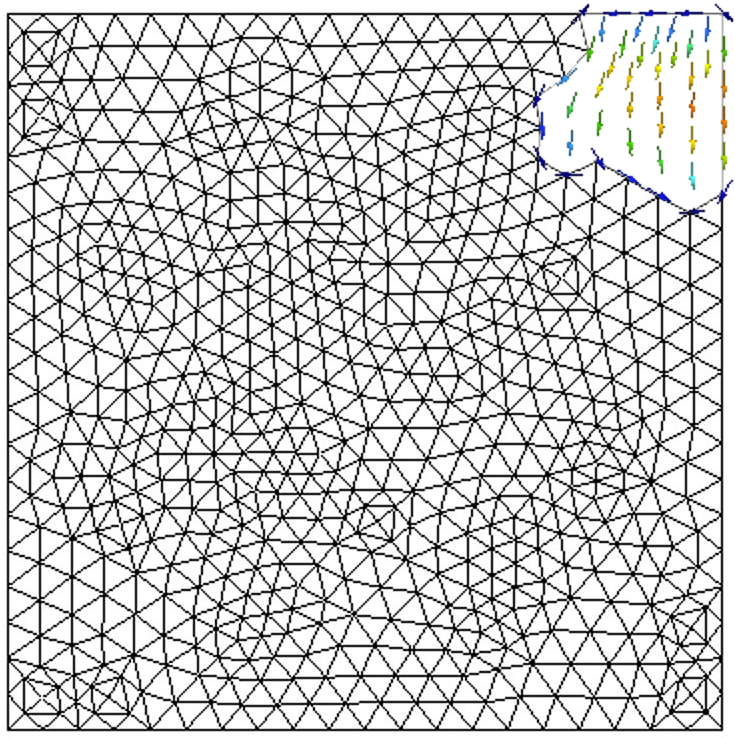}}\quad
\subfloat[][$L^2$ basis function]{\includegraphics[width=0.3\textwidth]{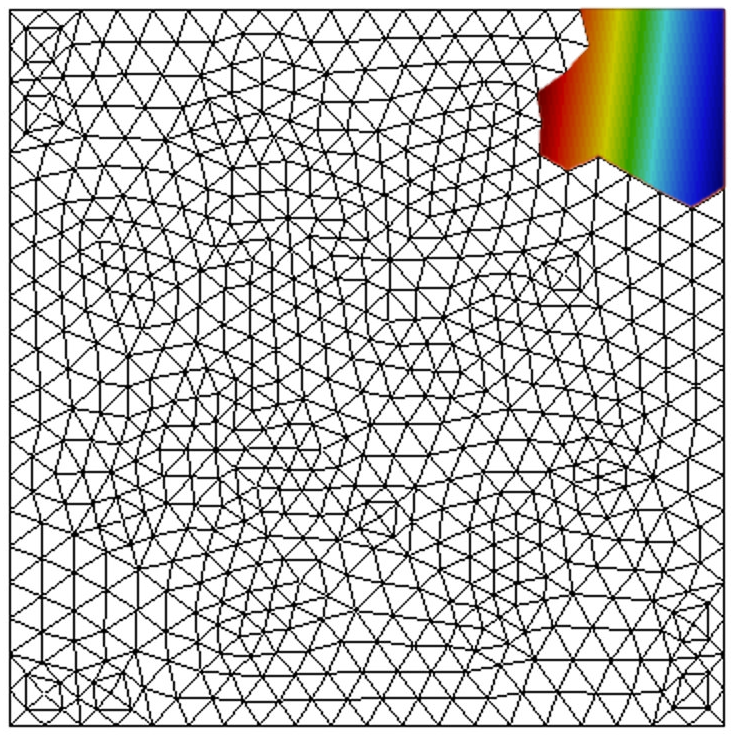}}
\caption[]{A two-dimensional illustration of coarse basis functions. Their supports match exactly the respective coarse elements.}
\label{fig:cbases}
\end{figure}

As indicated in \Cref{sssec:ebasedP}, the extension can be viewed in the local context of an agglomerate $T \in \cT^H$ and its associated lower-dimensional coarse entities. Note that some basis functions are supported on multiple agglomerates; see \Cref{fig:shapebasis}. Nevertheless, they are constructed by independent local agglomerate-by-agglomerate processes. The procedures executed on a single agglomerate $T$ produce the respective shape functions. Shape functions that form a particular basis function are conforming, i.e., they agree on fine-scale dofs shared between agglomerates. In the formation of the prolongation matrices, the final basis functions are obtained by joining together all associated shape functions. The shape functions on $T$ alone comprise the local-on-$T$ prolongation matrices, $P_{T,i}$; see \Cref{sssec:ebasedP}. Other basis functions are entirely supported in a single agglomerate $T$. As is customary, they are called \emph{bubble functions}; see \cite{BoffiMFE}. More specifically, for $i=1,\dots,3$, a $\cV^h(D_{T,i}\where T)$ bubble function is fully supported in $T$, and it is globally a function in $H(D_i\where \Omega)$, i.e., has a vanishing $\gamma_i$-trace (see \Cref{ssec:deRhamSpaces}) on $\partial T$. Here, $D_{X,i}$ denotes the restriction of the derivative operator $D_i$ ($i = 1,\ldots, 4$) to an entity $X \in \{V,E,F,T\}$, where $V$, $E$, $F$, and $T$ represent an agglomerated vertex, edge, facet, and element, respectively. \Cref{fig:cbases} illustrates coarse basis functions on agglomerates in two dimensions.

\subsubsection{Approximation property targets}
\label{sssec:ApproximationTargets}

The first step is to select so called \emph{targets}. The results in this paper are obtained using global polynomial targets due to their simplicity and inherent approximation properties following from standard polynomial approximation theory. Such polynomial targets are utilized to produce target traces on coarse entities (elements or agglomerates, facets, edges, vertices) and allow the construction of high-order (practically of arbitrary order) conforming and stable coarse spaces. The resulting target traces, obtained essentially by restricting the approximation property targets to the respective coarse entities, are then extended via the procedure in \Cref{sssec:extension}, providing the final coarse basis functions and coarse spaces with the desired order and approximation properties.

Global polynomial targets are set once in ParELAG on the finest level via simple call to the \verb|SetUpscalingTargets()| member function of the \verb|DeRhamSequenceFE| class. The procedure in ParELAG for building the targets is quite simple. On the finest level,  the finite element interpolants (via $\Pi^h_i$) of monomials up to a prescribed order constitute the targets. On coarse levels, the targets are transferred as needed via projection, i.e., by applying the operators $\Pi^H_i$, $i = 1,\dots,4$. Note that it is admissible to utilize polynomial targets of order higher than the order of the finest-level finite elements. In such a case, coarse basis functions have a high-order complexion but are represented piecewise by lower-order polynomials. Alternatively, local targets on coarse entities can be obtained by, e.g., solving local eigenvalue problems, which also provide approximation properties (cf. \cite{2016AMGeUpscaling,2012SAAMGe}). In any case, for the approach outlined here, appropriate respective target traces on coarse entities (elements or agglomerates, facets, edges, vertices) are obtained and available as needed, represented in terms of respective $\cV^h(D_i)$ dofs, $i=1,\dots,4$.

\subsubsection{PV traces and coarse trace spaces}
\label{sssec:PVtraces}

The lowest-dimensional traces for $\cV^H(D_i)$, $i=1,\dots,4$, from which extensions are initiated, are defined on agglomerated vertices, edges, facets, and elements, respectively. On these entities, the so called \emph{PV traces} (coming from \cite{2008deRhamAMGe}) are locally defined agglomerated entity by agglomerated entity. These traces alone provide, once extended following \Cref{sssec:extension}, the lowest-order stable coarse spaces on $\cT^H$.

Specifically, the PV traces 
$\phi_{\mathrm{PV}, 1}^V \in \cV^H(D_{V,1}\where V)$,
$\uphi_{\mathrm{PV}, 2}^E \in \cV^H(D_{E,2}\where E)$,
$\uphi_{\mathrm{PV}, 3}^F \in \cV^H(D_{F,3}\where F)$, and 
$\phi_{\mathrm{PV}, 4}^T \in \cV^H(D_{T,4}\where T)$ are defined to have unit integral on the corresponding agglomerated entity, $V$, $E$, $F$, and $T$. That is, the representation of the PV traces in terms of finer level dofs are given by
\begin{equation}
\label{eq:PV_target}
\begin{array}{ll}
\phi_{\mathrm{PV},1}^V = \phi_{\mathrm{PV},1}^v = 1 ,& v \equiv V \text{\small ($v$ is the fine vertex forming $V$)}\\
\uphi_{\mathrm{PV},2}^E\cdot\utau_E = \sum_{e\subset E} (\bvarphi^E)_e\;\uphi_{\mathrm{PV},2}^e\cdot\utau_e, & \text{satisfying } \int_E \, \uphi_{\mathrm{PV},2}^E\cdot\utau_E \,\dd \ell = 1, \\
\uphi_{\mathrm{PV},3}^F\cdot\un_F = \sum_{f\subset F} (\bvarphi^F)_f\;\uphi_{\mathrm{PV},3}^f\cdot\un_f, & \text{satisfying } \int_F \, \uphi_{\mathrm{PV},3}^F\cdot\un_F \,\dd \sigma = 1,  \\
\phi_{\mathrm{PV},4}^T = \sum_{\tau \subset T} \phi_{PV,4}^\tau, & \text{satisfying } \int_T \, \phi_{\mathrm{PV},4}^T \,\dd \Omega = 1,\\
\end{array}
\end{equation}
where $\utau_E$ and $\utau_e$ are the corresponding tangent vectors to the coarse edge $E$ and its constituting fine edges, $e$, $\un_F$ and $\un_f$ are the corresponding normal vectors to the coarse facet $F$ and its constituting fine facets, $f$, and the already available fine-scale PV traces $\phi_{\mathrm{PV},1}^v$, $\uphi_{\mathrm{PV},2}^e$, $\uphi_{\mathrm{PV},3}^f$, and $\phi_{PV,4}^\tau$ are utilized. Above, the quantities $(\bvarphi^E)_e$ and $(\bvarphi^F)_f$ denote the reciprocal orientation of the fine edge $e$ or fine facet $f$ and coarse edge $E$ or coarse facet $F$; see \Cref{sssec:Hmesh}. 

Next, the approximation targets are added to the coarse space after removing any linear dependence. Specifically, coarse basis functions orthogonal (with respect to the $L^2$ inner product on those entities) to the PV traces and among each other are generated using SVD and local mass matrices formulated on the respective agglomerated entities. Doing so, all coarse basis functions, apart from those stemming the PV traces, have a zero mean and the corresponding trace-space mass matrices produced via RAP on all algebraic levels are diagonal.

Note that this concludes the construction of the space $\cV^H(D_4) = \bigoplus_{T \in \mathcal{T}_H }\cV^H(D_{T,4}\where T)$ and no extension procedure is needed for this space.

\begin{figure}
\centering
\setlength\tabcolsep{0pt}
\newlength{\scaledtextwidth}
\setlength{\scaledtextwidth}{0.62683\textwidth}
\subfloat[][Extension of the PV trace on the agglomerated facet]{
\begin{tabular}{ccc}
\begin{minipage}{0.42\scaledtextwidth}\includegraphics[width=\textwidth]{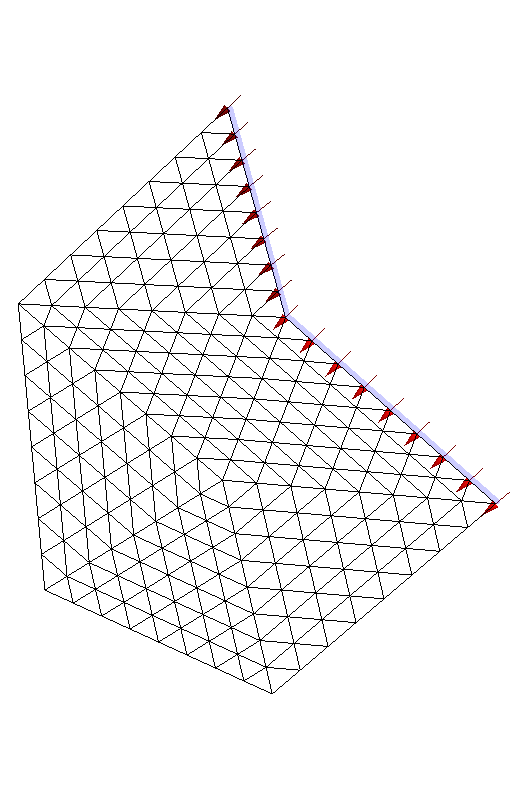}\end{minipage} &%
$\displaystyle\xrightarrow{\hphantom{.} \text{extension} \hphantom{.}}$ &%
\begin{minipage}{0.42\textwidth}\includegraphics[width=\textwidth]{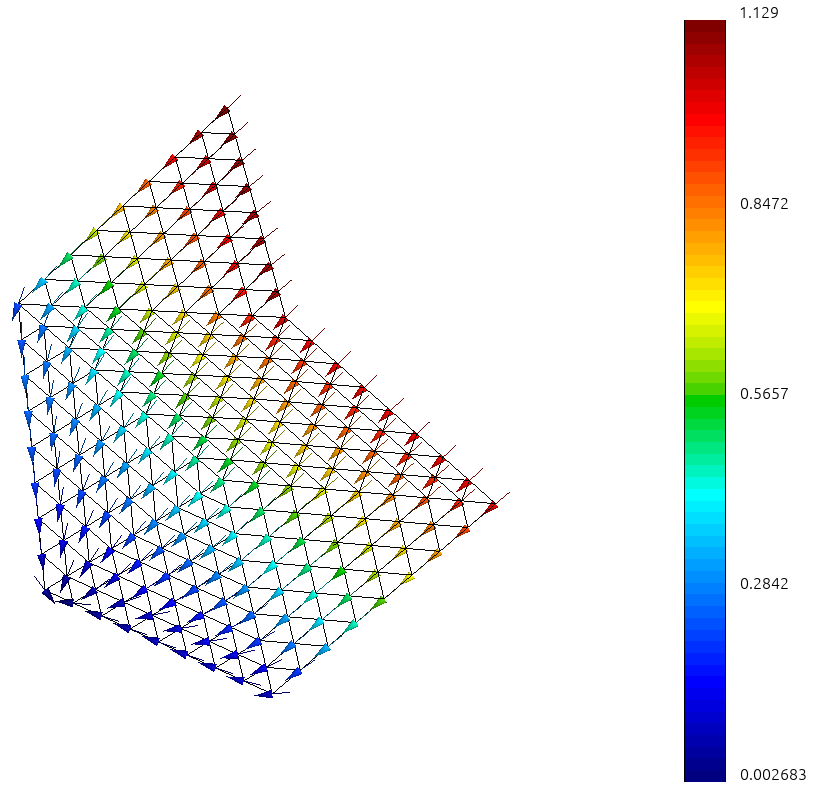}\end{minipage}
\end{tabular}\label{fig:extension_PV_facet}
}\\
\subfloat[][Extension of a trace, on the agglomerated facet, $L^2$-orthogonal to the PV trace]{
\begin{tabular}{ccc}
\begin{minipage}{0.42\scaledtextwidth}\includegraphics[width=\textwidth]{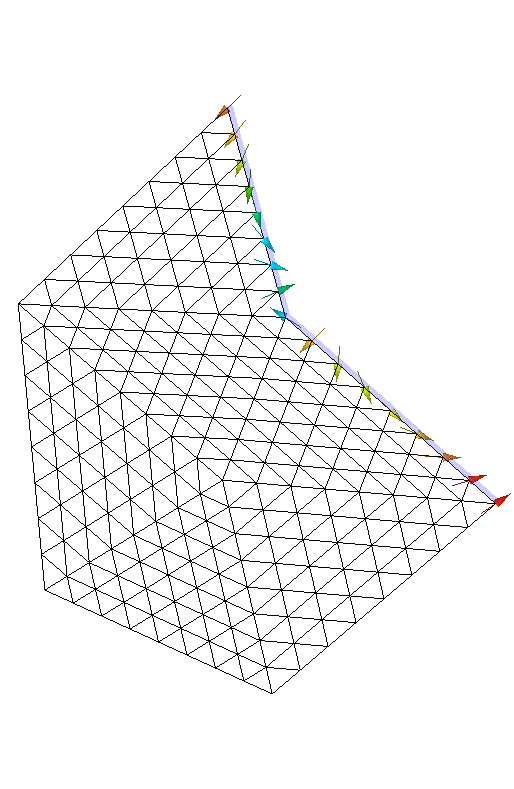}\end{minipage} &%
$\displaystyle\xrightarrow{\hphantom{.} \text{extension} \hphantom{.}}$ &%
\begin{minipage}{0.42\textwidth}\includegraphics[width=\textwidth]{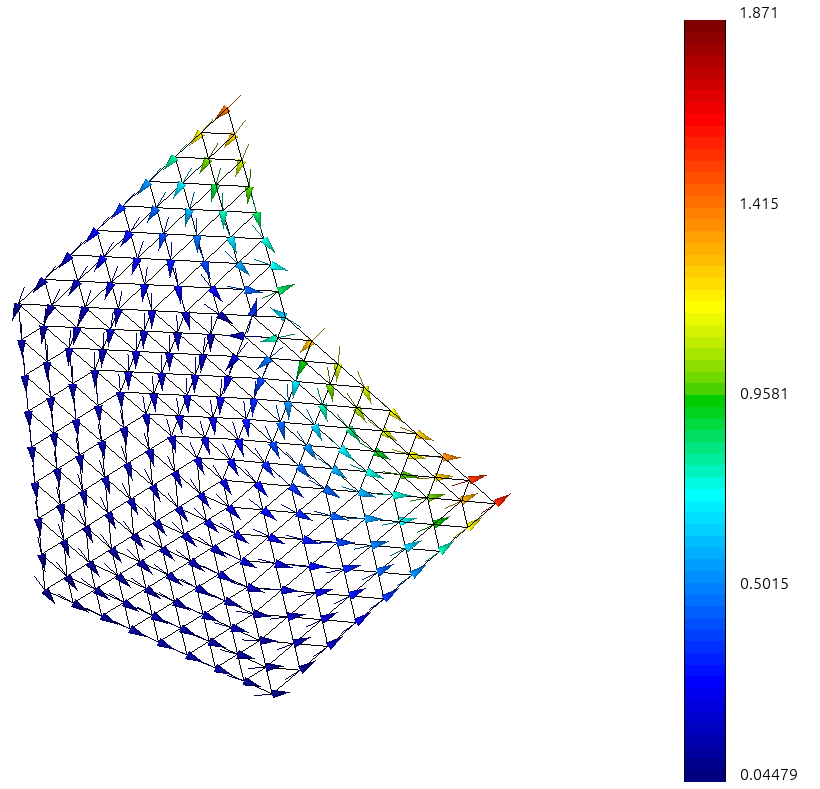}\end{minipage}
\end{tabular}\label{fig:extension_facet}
}\\
\subfloat[][Cross-space extension from an $L^2$ basis function to an $H(\div)$ bubble function]{
\begin{tabular}{ccc}
\begin{minipage}{0.42\textwidth}\includegraphics[width=\textwidth]{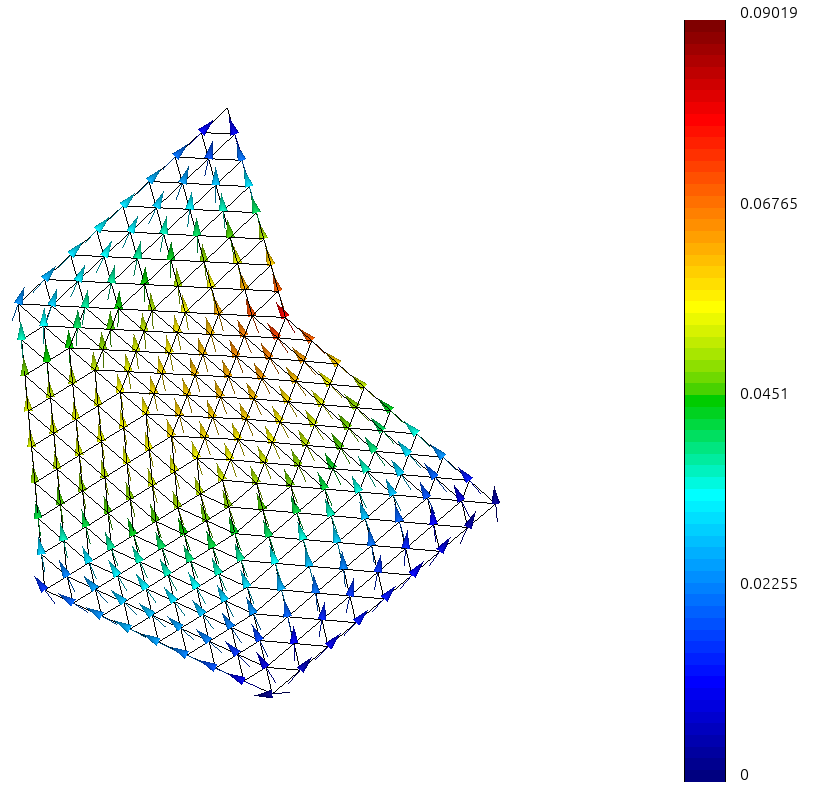}\end{minipage} &%
$\displaystyle\xleftarrow{\hphantom{.} \text{extension} \hphantom{.}}$ &%
\begin{minipage}{0.42\textwidth}\includegraphics[width=\textwidth]{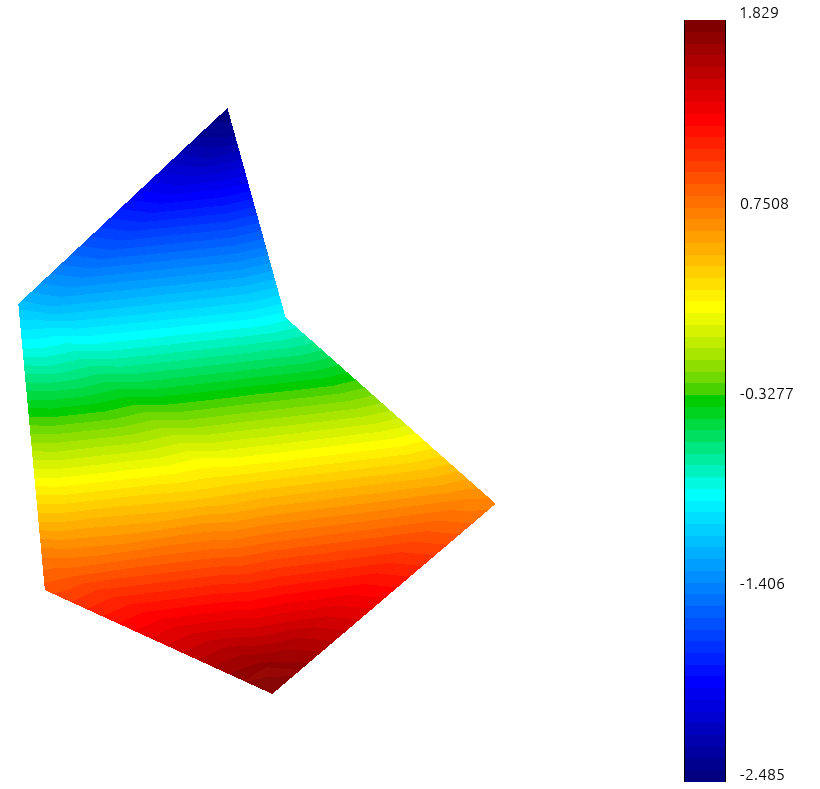}\end{minipage}
\end{tabular}\label{fig:extension_bubble}
}
\caption{A two-dimensional illustration of local extension procedures producing shape functions in $H(\div)$ on a sample agglomerated element, involving respective traces on an agglomerated facet (marked with a light shade) and bubble functions in the agglomerate.}\label{fig:extension}
\end{figure}

\subsubsection{Extension process}
\label{sssec:extension}

The extension procedure moves from right to left in the de Rham sequence \eqref{eq:HdeRham} and it considers two types of extensions: boundary extensions from lower to higher dimensional entities and cross-space extensions to ensure the exactness of the coarse sequence.  A detailed review of the extension process is provided here and an expository illustration is shown in \Cref{fig:extension}. For further analysis, including the feasibility of the extension problems, the demonstration of the exactness \eqref{eq:Hexact} and commutativity \eqref{eq:Hcomm} properties, see \cite{2014CoarseDeRham}.

\indent\paragraph{\textit{Extension from the lowest-dimensional traces}}

The first extension is from the lowest-dimensional, for the respective space, agglomerated entity to a one-dimension-higher agglomerated entity; see \Cref{fig:extension_PV_facet,fig:extension_facet} for an illustration. That is, for $i=1,\dots,3$, the extensions are respectively vertex to edge, edge to facet, and facet to element. The discussion here is associated with the method \verb|hFacetExtension()| of the class \verb|DeRhamSequence|, called within \verb|DeRhamSequence::Coarsen()|.

Let $L$ be a lowest-dimensional entity, $K$ a one-dimension-higher entity such that $L \subset \partial K$, and $\mu \in \cV^h(D_{L, i}\where L)$ be a given target trace on $L$. 
Denote with $D_{K,i}$ the restriction of the differential operator $D_i$ to the entity $K$, and with $D_{K,i}^*$ its adjoint.  The discretized version of $D_{K,i}$ is obtained by extracting the corresponding submatrix from $D_i^h$, while the discretized adjoint is obtained by matrix transposition. Also, introduce the trace operator $\gamma_{K,i}: \cV^h(D_{K,i}\where K) \to \cV^h(D_{\partial K,i}\where \partial K) $, which is defined in \Cref{ssec:deRhamSpaces} for an arbitrary domain $\Omega$. Namely, for $i=1,\dots,3$, $L$ is respectively an agglomerated vertex, edge, and facet, while $\gamma_{K,i}$ on $L$ is respectively a point value on the vertex, tangential flow on the edge, and normal flux on the facet.
Using the above notation, the extension of $\mu$ to $K$, $\phi_e \in \cV^h(D_{K,i}\where K)$, is obtained by solving the local PDE formally expressed as:
\begin{equation}
\label{eq:boundary_ext_1}
\left\{
\begin{array}{lll}
\phi_e + D_{K,i}^* \, \psi &= 0 &\text{ in } K,\\
D_{K,i} \, \phi_e &= \lambda \; \phi_{\mathrm{PV},i+1}^K \quad &\text{ in } K,\\
(\psi, \phi_{\mathrm{PV},i+1}^K)_{L^2(K)} & = 0 & \text{ in } \mathbb{R},\\
\gamma_{K,i} \, \phi_e &= \mu &\text{ on } L,\\
\gamma_{K,i} \, \phi_e &= 0 &\text{ on } \partial K \setminus L,
\end{array}
\right.
\end{equation}
where $\phi_{\mathrm{PV},i+1}^K$ is the PV trace in $\cV^h(D_{K, i+1}\where K)$ associated with $K$, $\psi$ is a local-on-$K$ function in $\cV^h(D_{K,i+1}\where K)$, and $\lambda$ is the scalar Lagrangian multiplier associated with the orthogonality constraint $(\psi, \phi_{\mathrm{PV},i+1}^K)_{L^2(K)} = 0$.
The orthogonality constraint guarantees the solvability of the above problem by ensuring that $\psi$ has zero mean on $K$. The value of $\lambda = (\mu, \phi_{\mathrm{PV},i}^L)_{L^2(L)}$ is determined by use of the Stokes' theorem thanks to \eqref{eq:PV_target}.

Next, $\cV^h(D_{K,i}\where K)$ bubble functions on $K$ are obtained by a cross-space extension from $\cV^h(D_{K,i+1}\where K)$ to $\cV^h(D_{K,i}\where K)$; see \Cref{fig:extension_bubble} for an illustration. This is necessary to preserve the exactness of the coarse sequence; see \eqref{eq:Hexact}. For each target trace $\phi_{\perp,i+1}^K\in \cV^h(D_{K,i+1}\where K)$ such that $(\phi_{\perp,i+1}^K, \phi_{{\mathrm{PV}},i+1}^K)_{L^2(K)} = 0$, the corresponding bubble function, $\phi_b \in \cV^h(D_{K,i}\where K)$, is obtained by solving
\begin{equation}
\label{eq:bubble_extension_1}
\left\{
\begin{array}{lll}
\phi_b + D_{K,i}^* \, \psi &= 0 &\text{ in } K,\\
D_{K,i} \, \phi_b &= \phi_{\perp,i+1}^K + c \; \phi_{\mathrm{PV},i+1}^K \quad &\text{ in } K,\\
(\psi, \phi_{\mathrm{PV},i+1}^K)_{L^2(K)} & = 0 & \text{ in } \mathbb{R}\\
\gamma_{K,i} \, \phi_b &= 0 &\text{ on } \partial K,
\end{array}
\right.
\end{equation}
where $\psi \in \cV^h(D_{K,i+1}\where K)$ is also zero-mean, $c = 0$, and the constraint serves to stabilize the system.

Finally, to ensure the approximation properties, additional $D_{K,i}$-free bubble functions\footnote{That is, functions $\phi \in \cV^h(D_{K,i}\where K)$ such that $D_{K,i} \phi = 0$ and $\gamma_{K,i} \phi = 0$.} in $\cV^h(D_{K,i}\where K)$ are produced by projecting the given target traces in $\cV^h(D_{K,i}\where K)$ associated with $K$ onto the respective space of $D_{K,i}$-free bubble functions and filtering out any linear dependence.

Upon completion of all the extensions in this steps, shape functions, co-chain projectors and exterior derivative operators for the spaces $\cV^H(D_{E,1}\where E)$, $\cV^H(D_{F,2}\where F)$, $\cV^H(D_{T,3}\where T)$ are defined for all coarse edges $E$, facets $F$, and elements $T$ of the coarse mesh. In particular, the construction of $\cV^H(D_3) = \bigoplus_{T \in \mathcal{T}_H} \cV^H(D_{T,3}\where T)$ is now complete.

\indent\paragraph{\textit{Further extensions to higher-dimensional agglomerated entities}}

For the case of $i=2$, one more extension step ($\text{facets}\to\text{elements}$) is necessary, while two such steps ($\text{edges}\to\text{facets}\to\text{elements}$) are required for $i=1$. Each such step has the form presented below. The discussion here, for each extension step, is associated with the method \verb|hRidgePeakExtension()| of the class \verb|DeRhamSequence|, called within \verb|DeRhamSequence::Coarsen()|.

Let $N$ be a lower-dimensional agglomerated entity (but not a lowest-dimensional one), $M$ a one-dimension-higher agglomerated entity such that $N \subset \partial M$, and $\eta \in \cV^h(D_{N,i}\where N)$ a trace on $N$ produced by a previous (lower-dimensional) extension. 
Let $s^\eta \in \cV^h(D_{M,i+1}\where M)$ be the extension of $D_{M,i} \eta$ from $N$ to $M$. Note the $s^\eta$ is known since the spaces $\cV^H(D_{N,i}  \where N)$ and $\cV^H(D_{M,i+1}  \where M)$ were already constructed during the previous extension step. Then, the respective extension of $\eta$ to $M$, $\phi_E$, in $\cV^h(D_{M,i}\where M)$ is obtained by solving the following formal local PDE:
\begin{equation}
\label{eq:boundary_ext_f}
\left\{
\begin{array}{lll}
\phi_E + D_{M,i}^* \, \chi &= 0 &\text{ in } M,\\
D_{M,i} \, \phi_E - D_{M,i+1}^* D_{M,i+1} \, \chi &= s^\eta \quad &\text{ in } M,\\
\gamma_{M,i} \, \phi_E &= \eta &\text{ on } N,\\
\gamma_{M,i} \, \phi_E &= 0 &\text{ on } \partial M \setminus N,
\end{array}
\right.
\end{equation}
where $\chi$ is a local-on-$M$ function in $\cV^h(D_{M,i+1}\where M)$. Above, the stabilization term $D_{M,i+1}^* D_{M,i+1} \, \chi$ is added to guarantee the uniqueness of  $\chi$.

To ensure the exactness \eqref{eq:Hexact}, $\cV^h(D_{M,i}\where M)$ bubble functions on $M$ must be included. These bubble functions, $\phi_B$, are such that $D_{M,i} \phi_B = \phi_{0,i+1}^M$, where $\phi_{0,i+1}^M \in \Ker(D^H_{M,i+1}) \subset \cV^H(D_{M,i+1}\where M)$ but expressed in practice in fine-scale $\cV^h(D_{M,i+1}\where M)$-dofs. Note that a basis of $\Ker(D^H_{M,i+1})$ was already constructed in the previous extension step and it is associated with respective $D_{M,i+1}$-free bubble functions and target traces with a zero mean (i.e., orthogonal to the respective PV targets) in $\cV^H(D_{M,i+1}\where M)$. Then, for each such $\phi_{0,i+1}^M$, the corresponding bubble function, $\phi_B$, is obtained by solving
\begin{equation}
\label{eq:bubble_ext_f}
\left\{
\begin{array}{lll}
\phi_B + D_{M,i}^* \, \chi &= 0 &\text{ in } M,\\
D_{M,i} \, \phi_B - D_{M,i+1}^* D_{M,i+1} \, \chi &= \phi_{0,i+1}^M \quad &\text{ in } M,\\
\gamma_{M,i} \, \phi_B &= 0 &\text{ on } \partial M.
\end{array}
\right.
\end{equation}
Note that the stabilization term $D_{M,i+1}^* D_{M,i+1} \, \chi$ is zero for this particular choice of right-hand side.

In the end, the given target traces in $\cV^h(D_i)$ associated with $M$ are projected on the space of $D_{M,i}$-free bubble functions on $M$ in $\cV^h(D_{M,i}\where M)$ and added towards the basis (possibly awaiting further extension) for $\cV^H(D_{M,i} \where M)$, after filtering out any linear dependence.

After one sweep of the above procedure, the space $\cV^H(D_2) = \bigoplus_{T\in \mathcal{T}_H} \cV^H(D_{T,2}\where T)$ is finalized, while an additional sweep is needed for constructing $\cV^H(D_1) = \bigoplus_{T\in \mathcal{T}_H} \cV^H(D_{T,1}\where T)$ and thus completing the coarse sequence. Finally,  note that approximation targets are not included in $\cV^H(D_1)$, since any such bubble function would vanish everywhere.

\subsubsection{Summary}
\label{sssec:summary}

The construction of the coarse shape functions for all spaces in the sequence is summarized in the diagram below, where $V$, $E$, $F$, and $T$ represent an agglomerated vertex, edge, facet, and element, respectively, and $D_{X,i}$ the restriction of the exterior derivative $D_i$ ($i = 1,\ldots, 4$) to an entity $X \in \{V,E,F,T\}$.
\begin{small}
\begin{equation*}
\begin{tikzcd}[column sep=2em]
\cV^H(D_{V,1} \where V) \arrow[d, "\eqref{eq:boundary_ext_1}"] \blacktriangleleft &                    & & \\
\cV^H(D_{E,1} \where E) \arrow[d, "\eqref{eq:boundary_ext_f}"] & \arrow[l, "\eqref{eq:bubble_extension_1}"]\cV^H(D_{E,2} \where E) \arrow[d, "\eqref{eq:boundary_ext_1}"] \circlearrowleft \blacktriangleleft &  & \\
\cV^H(D_{F,1} \where F) \arrow[d, "\eqref{eq:boundary_ext_f}"] & \arrow[l, "\eqref{eq:bubble_ext_f}"]\cV^H(D_{F,2} \where F) \arrow[d, "\eqref{eq:boundary_ext_f}"] \circlearrowleft & \arrow[l, "\eqref{eq:bubble_extension_1}"] \cV^H(D_{F,3}\where F) \arrow[d, "\eqref{eq:boundary_ext_1}"] \circlearrowleft\blacktriangleleft  & \\
\cV^H(D_{T,1} \where T) & \arrow[l, "\eqref{eq:bubble_ext_f}"] \cV^H(D_{T,2} \where T) \circlearrowleft & \arrow[l, "\eqref{eq:bubble_ext_f}"] \cV^H(D_{T,3}\where T) \circlearrowleft & \arrow[l, "\eqref{eq:bubble_extension_1}"] \cV^H(D_{T,4}\where T) \circlearrowleft \blacktriangleleft
\end{tikzcd}
\end{equation*}
\end{small}%
The symbols $\circlearrowleft$ and $\blacktriangleleft$ denote, respectively, the insertion---after filtering out linear dependencies---of the approximation property targets (\Cref{sssec:ApproximationTargets}) and the construction of the PV traces (\Cref{sssec:PVtraces}).
Vertical arrows denote boundary extensions from lower-dimensional entities to higher-dimensional ones, while horizontal arrows denote cross-space extensions necessary to ensure the exactness property \eqref{eq:Hexact}. The label next to each arrow represents the local PDE that is solved to compute the extension (cf. \Cref{sssec:extension}).

Given a finer level hierarchy, the construction of the coarse level sequence starts from the trace spaces on the main diagonal of the diagram and moves towards the bottom left corner of the diagram using the boundary and cross-space extension operators. Particularly, for any $i \in [1,4]$ and any sensible coarse entity $X$, coarse basis functions (or traces thereof) on $X$ are obtained, either from given targets or via local computational procedures like the extensions above, as vectors in $\cV^h(D_{X,i} \where X)$ on the respective finer dofs, thus building the coarse $\cV^H(D_{X,i} \where X)$ by producing its basis and the corresponding coarse $\cV^H(D_{X,i} \where X)$-dofs associated with $X$.

\begin{remark}
In ParELAG, the extension operators in \eqref{eq:boundary_ext_1}, \eqref{eq:bubble_extension_1}, \eqref{eq:boundary_ext_f}, \eqref{eq:bubble_ext_f} can be attuned to the particular problem to be solved. For example, the coefficients $\alpha$ and $\beta$ in \eqref{eq:bfs} are incorporated in the appropriate local PDEs. Thus, the extension procedure can be informed about the particular problem of interest and the obtained coarse bases become problem-dependent.
\end{remark}

\section{Smoothers, coarse solvers, and the multigrid}
\label{sec:smoothsolve}

\begin{algorithm}
\caption{A procedure implementing a single multilevel V-cycle. It computes the action of a multilevel preconditioner $B^\mone_{\mathrm{ML}}$, i.e., $\bx_{\mathrm{ML}} = \bx_0 + B^\mone_{\mathrm{ML}}(\bb - A\bx_0)$, applied for solving the linear system $A\bx = \bb$, formulated on the finest level. The hierarchy of system matrices, $\{A_k\}_{k=1}^{\ell-1}$, where $A_1 = A$ denotes the fine-level system matrix, is precomputed before the invocation of the V-cycle. In particular, the coarser-level system matrices $A_k = (P^{k-1}_{k})^T A_{k-1}\, P^{k-1}_{k}$ for $k=2,\dots,\ell-1$ are obtained via a Galerkin RAP procedure.}\label{alg:ML}
\begin{algorithmic}
\small
\STATE \textbf{PROCEDURE}: $\bx_{\mathrm{ML}} \leftarrow \mathbf{ML}\lp \{A_k\}_{k=1}^{\ell-1},\; \bb,\; \bx_0,\; \{M_k\}_{k=1}^{\ell-1},\; \{P_{k+1}^k\}_{k=1}^{\ell-1},\; B^\mone_\ell,\; l \rp$
\STATE \textbf{INPUT}: A hierarchy of linear-system matrices $\{A_k\}_{k=1}^{\ell-1}$, a right-hand side vector $\bb$, a current iterate $\bx_0$, a hierarchy of relaxation $\{M_k\}_{k=1}^{\ell-1}$ and prolongation $\{P_{k+1}^k\}_{k=1}^{\ell-1}$ operators, a solver or preconditioner $B^\mone_\ell$ on the coarsest level, and a current level $l$.
\STATE \textbf{OUTPUT}: A new multigrid iterate $\bx_{\mathrm{ML}} \leftarrow \bx$.
\STATE \textbf{STEPS}:
\STATE Initialize: $\bx \leftarrow \bx_0$. 
\STATE Pre-relax: $\bx \leftarrow \bx + M_l^{-1} \lp \bb - A_l\bx \rp$.
\STATE Correct (evoke $B^\mone_\ell$ or recur):
\begingroup
\setlength{\itemindent}{1em}
\addtolength{\algorithmicindent}{1em}
\IF{$l = \ell - 1$ (i.e., coarsest level reached)}
\STATE $\be_c \leftarrow B^\mone_\ell (P^l_{l+1})^T (\bb - A_l\bx)$;
\ELSE
\STATE $\be_c \leftarrow \mathbf{ML}\lp \{A_k\}_{k=1}^{\ell-1},\; (P^l_{l+1})^T (\bb - A_l\bx),\; \bzero,\; \{M_k\}_{k=1}^{\ell-1},\; \{P_{k+1}^k\}_{k=1}^{\ell-1},\; B^\mone_\ell,\; l+1 \rp$;
\ENDIF
\STATE $\bx \leftarrow \bx + P^l_{l+1}\, \be_c$.
\endgroup
\STATE Post-relax: $\bx \leftarrow \bx + M_l^{-T} \lp \bb - A_l\bx \rp$.
\end{algorithmic}
\end{algorithm}

This section is devoted to describing a scalable multigrid preconditioner for the finite element matrices stemming from discretizations of the bilinear forms in \eqref{eq:bfs} attuned to the AMGe hierarchies of de Rham sequences as constructed by ParELAG (see \Cref{sec:HdeRham}).

Multigrid preconditioners, implemented via multilevel cycles such as the well-known V-cycle in \Cref{alg:ML} (see \cite{VassilevskiMG}), have three main components: a hierarchy of spaces given in the form of a hierarchy of prolongator/restriction operators, a \emph{relaxation} (or \emph{smoothing}) procedure, and a solver or preconditioner for the coarsest problem. 

In \Cref{alg:ML}, $\ell$ denotes the number of levels in the hierarchy, and $l$ the current level in the hierarchy. As usual in algebraic multigrid literature, the finest level (where the matrix $A$ and right-hand side $\bb$ are defined) is denoted by $l=1$. The smoothing procedure and prolongator operator at level $l$ in the hierarchy are denoted by $M_l$ and $P_{l+1}^l$, respectively, while the coarse grid solver is denoted by $B_\ell^{\mone}$. Externally, the procedure is invoked with $l=1$ and it calls itself recursively for an increasing value of $l$, until the coarsest level $l=\ell$ is reached. Note that while in principle \Cref{alg:ML} can be used standalone iteratively to obtain a stationary (or fixed-point) iterative method, the main interest here is to apply the preconditioner $B_{\mathrm{ML}}^\mone$ within a preconditioned conjugate gradient (PCG) method. In that case, \Cref{alg:ML} is invoked at each PCG iteration with $\bx_0 = \bzero$ and $\bb$ denoting the residual at the current iteration.

In what follows, $A_{D_i}$ denotes the finite element matrix stemming from a discretization of the bilinear form $a_{D_i}(u,v) = (\alpha_i\, D_i u, D_i v)_0 +(\beta_i\, u, v)_0$ for $u,v \in H(D_i)$ and $\alpha_i > 0$, $\beta_i \ge 0$. The cases $i=2,3$ corresponds, respectively, to the $H(\ucurl)$ and $H(\div)$ forms in \eqref{eq:bfs} considered here. By use of the ParELAG prologation operators $\{P_{k+1}^k\}_{k=1}^{l-1}$, the Galerkin projections of $A_{D_i}$ at a generic coarse level $l+1$ are defined as $A_{D_i}^{l+1} = (P_{l+1}^l)^T A_{D_i}^{l} P_{l+1}^l$, where $A_{D_i}^{1} = A_{D_i} = A_{D_i}^h$ corresponds to the finest level (i.e., the finite element level).

The constructions of hybrid (``combined,'' a.k.a. ``Hiptmair'') smoothers $M_l$ for $A_{D_i}^{l}$ and coarse AMS/ADS solvers $B_\ell^{-1}$ for $A_{D_i}^{\ell}$ are respectively described in \Cref{ssec:hsmoothers,ssec:AMSADS}. For the coarse solvers, in particular, ParELAG constructs the transition and projection operators on the AMGe coarse de Rham sequence, which mimic those of the fine-grid de Rham sequence constructed in MFEM, for use within the AMS and ADS solvers in the HYPRE library \cite{hypre}. 

\subsection{Relaxation via hybrid smoothers}
\label{ssec:hsmoothers}

A general level-independent smoothing procedure based on \cite{1998Hiptmair} (see also \cite[Appendix F]{VassilevskiMG}) is outlined here, providing the relaxation processes $\{M_k\}$ in \Cref{alg:ML} for all levels.

For each matrix $A_{D_i}^l$, let $M_{D_i}^l$ denote the corresponding (point) smoother, e.g., a Jacobi-type, Gauss-Seidel-type, or their block or $\ell^1$-scaled variants smoothers. It is well known \cite{2000HdivHcurlMultigrid} that these type of smoothers do not lead to optimal (mesh independent) multigrid preconditioners for $H(\ucurl)$ and $H(\div)$ forms, due to the large near-null-space of these operators. Hybrid smoothers need to ``reach'' in the reverse direction of the de Rham sequence to perform an additional smoothing step on the near-null-space components.

Specifically, hybrid smoothers leverage certain decompositions of the spaces and the exactness of the de Rham sequence \eqref{eq:Hexact}, similar to the methods in \Cref{ssec:AMSADS} below. For $i=2,3$, at a generic level $l$ of the hierarchy, the stable decomposition $\cV^l(D_i) = \Ker(D_i^l) \oplus [\Ker(D_i^l)]^\perp$ allows to define efficient smoothers for each component separately. Smoothing the $[\Ker(D_i^l)]^\perp$ component is addressed by $M_{D_i}^l$. Smoothing the component in $\Ker(D_i^l)$ requires the construction of a point smoother $M_{D_{i-1}}^l$ for the \emph{auxiliary} matrix $(D_{i-1}^l)^T A_{D_i}^l D_{i-1}^l$, where the transition operator $D_{i-1}^l$ is an exterior derivative at level $l$ constructed in ParELAG.

Then, for $i=2,3$, the smoothing step $\bx_1 = \bx_0 + (\bbM_{D_i}^l)^\mone (\bb - A_{D_i}^l\bx_0)$ with the combined smoother $(\bbM_{D_i}^l)$ reads
\begin{equation}\label{eq:hsmoother}
\begin{split}
\bx_\half &= \bx_0 + (M_{D_i}^l)^\mone(\bb - A_{D_i}^l\bx_0),\\
\bx_1 &= \bx_\half + D_{i-1}^l\, (M_{D_{i-1}}^l)^\mone\, (D_{i-1}^l)^T (\bb - A_{D_i}^l\bx_\half).
\end{split}
\end{equation}

That is, the \emph{error propagation} operator satisfies
\[
I - (\bbM_{D_i}^l)^\mone A_{D_i}^l = \left [I - D_{i-1}^l\, (M_{D_{i-1}}^l)^\mone\, (D_{i-1}^l)^T\, A_{D_i}^l \right] \left [I - (M_{D_i}^l)^\mone A_{D_i}^l \right ].
\]
Notice that, since $D_i^l D_{i-1}^l = 0$ (see \eqref{eq:Hexact}), recursively utilizing $\bbM_{D_{i-1}}^l$ in place of $M_{D_{i-1}}^l$ in \eqref{eq:hsmoother} changes nothing. Therefore, there is no need to ``reach'' further than one step backwards into the de Rham sequence. To compute an iteration with $(\bbM_{D_i}^l)^{-T}$, reverse the order of the steps in \eqref{eq:hsmoother}, while respectively using $(M_{D_i}^l)^{-T}$ and $(M_{D_{i-1}}^l)^{-T}$ in place of the ones in \eqref{eq:hsmoother}.

In the numerical results presented in \Cref{sec:numerical}, the so called $\ell^1$-scaled symmetric block Gauss-Seidel smoother \cite{2011Smoothers}, as implemented in HYPRE, is used for all $M_{D_i}^l$. For more details on the analysis of the hybrid approach in a multigrid setting, which counts on the exactness property \eqref{eq:Hexact}, see \cite[Appendix F]{VassilevskiMG}. 

\subsection{Coarse solvers using AMS and ADS}
\label{ssec:AMSADS}

Similarly to \Cref{ssec:hsmoothers}, special (a.k.a. regular) decompositions (cf. \cite{2007AuxiliarySpace}) of the finite element spaces of interest are used to break the problem of obtaining a holistic preconditioner into the preconditioning of each component of the decomposition. This provides an \emph{auxiliary space preconditioner} that reduces the problem to preconditioning a few $H^1$-type forms, which can be efficiently addressed by AMG, and smoothing. As a part of HYPRE, BoomerAMG is used in this case. In this work, AMS and ADS, possibly wrapped in PCG and performing multiple iterations up to a given tolerance, are to be used as coarse solvers $B_\ell^\mone$ in \Cref{alg:ML}. 

To understand how to generalize AMS and ADS to the coarse problems generated by ParELAG, a brief overview of these methods applied to the finite element discretization level is presented below.
To this aim, the finite element space $\ucV^h(D_1) = \ucV^h(\ugrad) = [\cV^h(\ugrad)]^3$ of vectorial $H^1$-conforming functions is introduced. Next, for $i=2,3$ the interpolation operators $
\widehat{\Pi}_i^h: \ucV^h(D_1) \to \cV^h(D_i)$ are constructed by means of an overwriting assembly procedure collecting the local (element-by-element) versions of these operators.

Using the exterior derivative $D_{i-1}^h$ and the newly introduced interpolator $\widehat{\Pi}_i^h$, an arbitrary function $\bv^h \in \cV^h(D_i)$ ($i=2,3$) admits the stable decomposition \cite{2007AuxiliarySpace}
\begin{equation*}
\bv^h = \btv^h + \whPi^h_i \br^h + D_{i-1}^h \bz^h\quad\text{ for } \bv^h \in \cV^h(D_i), i=2,3,
\end{equation*}
where $\btv^h \in \cV^h(D_i)$, $\br^h \in \ucV^h(\ugrad)$, and $\bz^h \in \cV^h(D_{i-1})$.
Based on the above decomposition, (an additive version of\footnote{Several different variations, including multiplicative and ones that treat $\ucV^h(\ugrad)$ in a scalar component-wise fashion, are implemented in HYPRE \cite{hypre}.}) the auxiliary space preconditioner has the form \cite{2009AMS, 2012ADS}
\begin{equation}\label{eq:Bcurldiv}
(B_{D_i}^h)^\mone = (M_{D_i}^h)^\mone + \whPi^h_i\, (B_{H^1}^h)^\mone\, (\whPi^h_i)^T + D_{i-1}^h\, (B_{D_{i-1}}^h)^\mone\,(D_{i-1}^h)^T,
\end{equation}
where $M_{D_i}^h$ is a smoother for $A_{D_i}^h$ and $B_{H^1}^h$ is an AMG preconditioner (e.g., BoomerAMG) of the vector $H^1$-type matrix $A_{H^1}^h = (\whPi^h_i)^T A_{D_i}^h \whPi^h_i$, and $B_{D_{i-1}}^h$ is a multilevel preconditioner for $A_{D_{i-1}}^h = (D_{i-1}^h)^T A_{D_i}^h  D_{i-1}^h$. Note that, for $i=2$, the matrix $A_{D_{i-1}}^h$ is equivalent to the discretization of an $H^1$-conforming form, and thus can be preconditioned using AMG. For $i=3$, $A_{D_{i-1}}^h$ is equivalent to the discretization of a singular $H(\ucurl)$-conforming form, thus requiring the use of another auxiliary space preconditioner  \eqref{eq:Bcurldiv} with $i=2$.

The auxiliary space preconditioner \eqref{eq:Bcurldiv} can be seamlessly generalized to the preconditioning the coarse matrices $A_{D_i}^\ell$ ($i=2,3$), generated using the ParELAG hierarchy of nested de Rham sequences, by replacing the finite element exterior derivative $D_{i-1}^h$ and interpolation operator $\widehat{\Pi}_i^h$ with their coarse level counterparts. In particular, the coarse level interpolation operator $\whPi_i^\ell: \ucV^\ell(\ugrad) \to \cV^\ell(D_i)$ is constructed in ParELAG by means of a RAP procedure from the finer level interpolator. That is, denoting the vectorial counterpart of the ParELAG prolongator $(P_1)_{l+1}^{l}: \cV^{l+1}(D_1) \to \cV^{l}(D_1)$ by $\widehat{P}_{l+1}^l: \ucV^{l+1}(D_1) \to \ucV^{l}(D_1)$, $\whPi_i^\ell$ is obtained by applying $\ell-1$ times the recursion
\begin{equation*}
    \whPi_i^{l+1} = (\Pi_i)_{l+1}^l\, \whPi_i^l\, \widehat{P}_{l+1}^l, \quad l=1, \ldots, \ell-1,
\end{equation*}
where $\whPi_i^1 = \whPi_i^h$ is the interpolation operator at the finite element discretization level (provided by MFEM) and $(\Pi_i)_{l+1}^l$ denotes the respective ParELAG-generated co-chain projector from level $l$ onto level $l+1$.

\subsection{Overview of composite solvers in ParELAG}

The ParELAG library provides access to a variety of solvers for sparse linear systems, including for block systems. Some of them are implemented within ParELAG itself, like the hybrid smoothers of \Cref{ssec:hsmoothers} and the V-cycle of \Cref{alg:ML} for AMGe, while others make use of external libraries, such as HYPRE \cite{hypre}, MFEM \cite{mfem}, SUPERLU\_DIST \cite{superludist}, and STRUMPACK \cite{strumpack}. Furthermore, solvers can be combined into a composite solver for the linear system of interest.

ParELAG achieves this by first generating a \emph{solver (or preconditioner) library} (an object of class \verb|SolverLibrary|) from an XML configuration file with a very intuitive syntax. Such a file declares all the solvers and preconditioners needed by the application, together with their specific parameters and how they are combined. A solver is declared by assigning a name and a list of parameters for a particular method internally provided by ParELAG. For example, one can declare a solver in the XML file that represents a multigrid method like the one in \Cref{alg:ML} and appoint other solvers from the solver library to act as smoothers and coarse solvers, which have their own sets of parameters and may, in turn, internally employ other solvers or preconditioners from the solver library. The static member function \verb|CreateLibrary()| is defined to instantiate a \verb|SolverLibrary| object from the provided XML configuration file.

To instantiate a solver, users first interrogate the \verb|SolverLibrary| by calling its member function \verb|GetSolverFactory()|, which returns a \emph{solver factory} (an object of class \verb|SolverFactory|) for the desired solver. The solver (i.e., an object of class \verb|Solver|\footnote{\verb|Solver| is a virtual class defined in MFEM.}) is then instantiated by calling the member function \verb|BuildSolver()| of \verb|SolverFactory|. For solving a linear system, users then call the \verb|Mult()| member, as they would do with any linear solver implemented in MFEM. Such a paradigm of solver structuring may sound familiar to a reader that has been exposed to some of the popular available solver libraries.

Examples of XML configuration files for solving the $H(\ucurl)$ and $H(\div)$ problems presented in \Cref{sec:numerical} are included in the ParELAG mini application of MFEM. Additional examples to configuring solver factories for AMS, ADS, Krylov space methods, hybrid smoothers, AMGe cycles, block preconditioners, and other methods within HYPRE and MFEM, can be found in the ParELAG library.

\section{Numerical examples}
\label{sec:numerical}

This section contains numerical results employing ParELAG in the context of the discussed AMGe multigrid solvers for \eqref{eq:bfs}. The results are produced using the ParELAG mini application \cite{miniapps} in MFEM.

\subsection{On the benchmark problem}

\begin{figure}
\centering
\captionsetup[subfloat]{labelformat=empty,width=0.25\textwidth}
\subfloat{\includegraphics[width=0.25\textwidth]{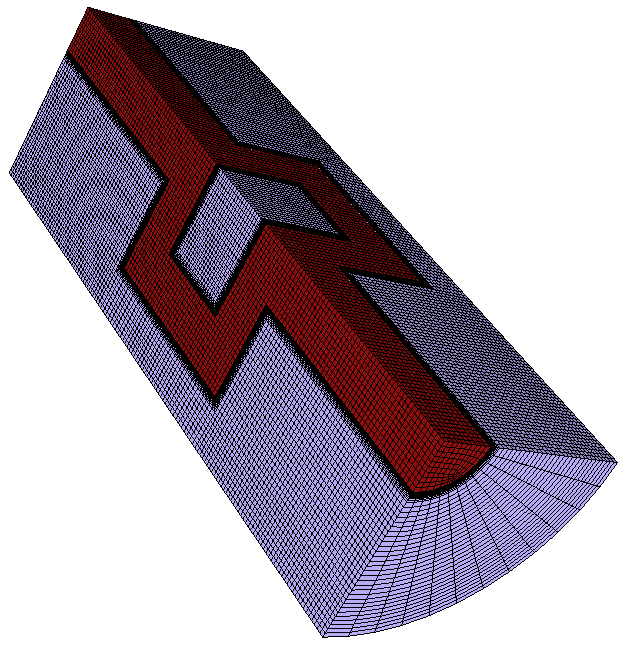}}\quad\quad\quad\quad
\subfloat{\includegraphics[width=0.25\textwidth]{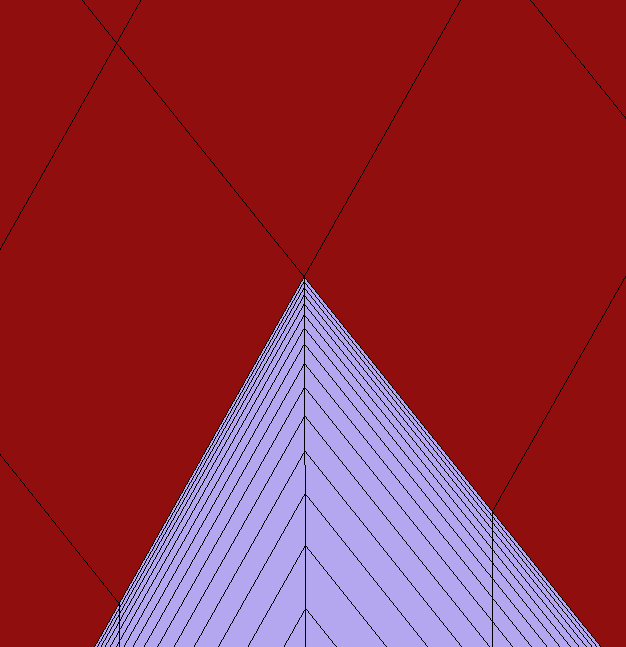}}\\
\subfloat[][$(\alpha, \beta) = (1.641, 0.2)$]{\includegraphics[width=0.25\textwidth]{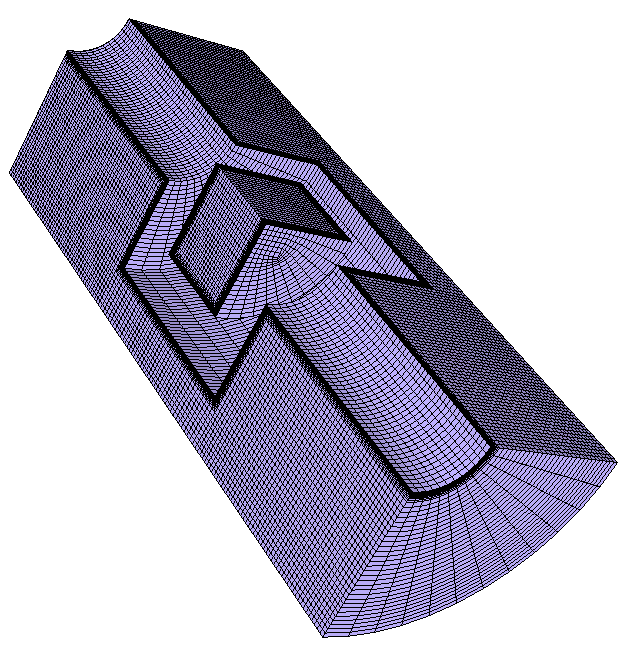}}\quad\quad\quad\quad
\subfloat[][$(\alpha, \beta) = (0.00188, 2000)$]{\includegraphics[width=0.25\textwidth]{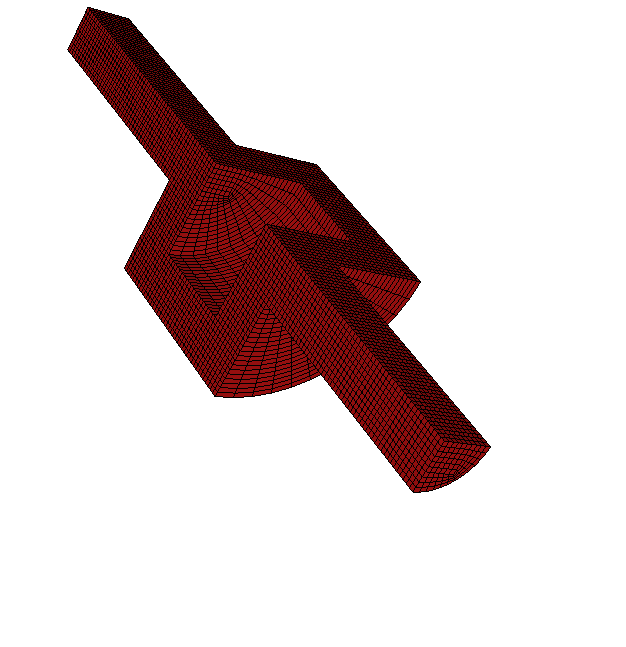}}
\caption[]{Domain, initial or starting mesh (including a close-up of the graded mesh in the upper right) of 116640 hexahedral elements, and piecewise constant coefficients for the numerical examples.}\label{fig:coeff}
\end{figure}

The benchmark considered here is inspired by the so-called ``crooked pipe'' problem (see \cite{graziani2000crooked, gentile2001implicit, 2012ADS}). This involves solving $H(\ucurl)$ and $H(\div)$ forms with scalar discontinuous coefficients with large jumps on a graded mesh with highly stretched (anisotropic) elements. For demonstration, the methodology is applied for solving linear systems coming from discretizations of formulations using the bilinear forms in \eqref{eq:bfs}, a constant right-hand side, and homogeneous essential boundary conditions for simplicity. The computational domain comprised of two different materials and the (coarsest) finite element mesh are depicted in \Cref{fig:coeff}. The coefficients are set as $(\alpha, \beta) = (1.641, 0.2)$ in the outer material (depicted with the lighter color) and $(\alpha, \beta) = (0.00188, 2000)$ in the inner core (depicted with the darker color). 

The systems coming from \eqref{eq:bfs} are solved using inexact PCG (see, e.g., \cite{golub1999inexact}) preconditioned by a single AMGe V-cycle; see  \Cref{sec:smoothsolve}. The iterative process is stopped when the relative size of the residual, measured by the preconditioner-induced norm, is reduced by six order of magnitude (i.e., $10^{-6}$ relative tolerance).

The V-cycles use AMGe hierarchies, as described above, and hybrid smoothers for pre and post-relaxation; see \Cref{ssec:hsmoothers}. An application of the hybrid smoother invokes two sweeps of $\ell^1$-scaled symmetric block Gauss-Seidel for each \emph{primary} and \emph{auxiliary} smoothings within the hybrid approach. A fixed number of five iterations of PCG preconditioned by AMS or ADS, respectively, serves as a solver on the coarsest level\footnote{The particular choice is largely motivated by the objective to demonstrate the flexibility of the ParELAG \verb|SolverLibrary| and its ability to combine a variety of solvers and smoothers. Using a single or a few applications of the respective AMS or ADS, without PCG, is also a valid option here.}; see \Cref{ssec:AMSADS}.

In the tests, the mesh in \Cref{fig:coeff} is uniformly refined multiple times to obtain a fine-grid problem, which is consequently solved in parallel by the methods discussed in this paper. \emph{Weak scaling} experiments are considered here. That is, as the mesh is refined the number of processors is also increased so that the number of elements per processor is maintained constant.

\subsection{Results}

Computational results on solving systems induced by \eqref{eq:bfs} in parallel for these generally challenging problems are presented here, employing lowest and next to the lowest order finite elements and uniformly refining the initial mesh in \Cref{fig:coeff}. The tests utilize the Quartz cluster at Lawrence Livermore National Laboratory. It is equipped on each node with 128 GB of memory and two 18-core Intel Xeon E5-2695 v4 CPUs at 2.1 GHz, resulting in 36 computational cores per node, and the total number of computational nodes (cores) is 2,988 (107,568). The peak single CPU memory bandwidth is 77 GB/s and the Cornelis Networks Omni-Path provides the inter-node connection.

The number of PCG iterations ($\mathrm{it}_\mathrm{e}$), the number of dofs, and the number of elements (elems) on the finest level are reported, as well as the number of uniform mesh refinements (refs) employed to obtain the fine mesh, the total number of levels (denoted by $\ell$) in the AMGe hierarchy, and the number of utilized processors\footnote{Strictly speaking, this is the number of individual independent computational units, i.e., cores.} (procs). Also, the \emph{grid complexity} (GC) is reported, which is the total number of dofs in the hierarchy, respectively of nested subspaces of $H(\ucurl)$ or $H(\div)$, over the number of finest dofs. The \emph{operator complexity} (OC) is defined similarly but using the number of nonzeros in the sparsity patterns of the matrices in the AMGe hierarchy instead.

For comparison, results obtained with PCG preconditioned by HYPRE AMS and ADS acting on the finite element level are also presented. These includes the number of iterations ($\mathrm{it}_\mathrm{a}$) and wall-clock timings. Also, for comparison and completeness, two test cases are demonstrated: one where the number of AMGe levels is kept fixed (equal to 3) as the fine mesh is refined, essentially also refining the coarsest level, and another where as the fine mesh is refined the number of levels is increased so that the coarsest level is constant and coinciding with the initial mesh in \Cref{fig:coeff}.

Recall that the construction of the whole coarse de Rham sequences includes the element agglomeration times, the local extension procedures, and building other necessary constructs. As expected, we observe that the construction time almost does not grow when refining the mesh and increasing the number of processors to maintain constant number of elements per processor, especially in the case of constant number of levels, since the majority of the time spent is on the local extension procedures, which scale perfectly (they are embarrassingly parallel), as they involve no communication. Therefore, the weak scaling of the construction of the whole coarse de Rham sequences on all levels is almost perfect, i.e., it takes almost constant time.

\subsubsection{Results for lowest order elements}

\begin{table}
\centering
\footnotesize
\begin{tabular}{ | c | r | r | c | r | r | r | r | r | r | }
\hline
\multirow{2}{*}{refs} & \multirow{2}{*}{procs} & \multicolumn{1}{|c|}{\multirow{2}{*}{elems}} & \multirow{2}{*}{$\ell$} & \multicolumn{3}{|c|}{$H(\ucurl)$} & \multicolumn{3}{|c|}{$H(\div)$} \\
\cline{5-10}
& & & & \multicolumn{1}{|c|}{dofs} & $\mathrm{it}_\mathrm{e}$ & $\mathrm{it}_\mathrm{a}$ & \multicolumn{1}{|c|}{dofs} & $\mathrm{it}_\mathrm{e}$ & $\mathrm{it}_\mathrm{a}$ \\
\hline
2 & 72 & 7,464,960 & 3 & 22,772,484 & 131 & 93 & 22,583,232 & 36 & 31 \\
\hline
3 & 576 & 59,719,680 & 3 & 180,667,656 & 198 & 118 & 179,912,448 & 55 & 44 \\
\hline
4 & 4,608 & 477,757,440 & 3 & 1,439,303,184 & 231 & 148 & 1,436,285,952 & 86 & 60 \\
\hline
3 & 576 & 59,719,680 & 4 & 180,667,656 & 169 & 118 & 179,912,448 & 54 & 44 \\
\hline
4 & 4,608 & 477,757,440 & 5 & 1,439,303,184 & 190 & 148 & 1,436,285,952 & 77 & 60 \\
\hline
\end{tabular}
\caption[]{Weak scaling results with lowest order finite elements for both $H(\ucurl)$ and $H(\div)$ problems, as provided by \eqref{eq:bfs} and \Cref{fig:coeff}. Here, $\ell$ denotes the number of levels in the AMGe hierarchy, $\mathrm{it}_\mathrm{e}$ -- the number of PCG iteration using the proposed AMGe preconditioner, and $\mathrm{it}_\mathrm{a}$ -- the number of PCG iteration using the auxiliary space (AMS or ADS) preconditioners from the HYPRE library. In all cases, elems / procs = 103,680, while the GC and OC of the AMGe hierarchy are approximately constant and equal to 1.14 for all refinement levels.}\label{tbl:lo}
\end{table}

\begin{figure}
\centering
\includegraphics[width=\textwidth]{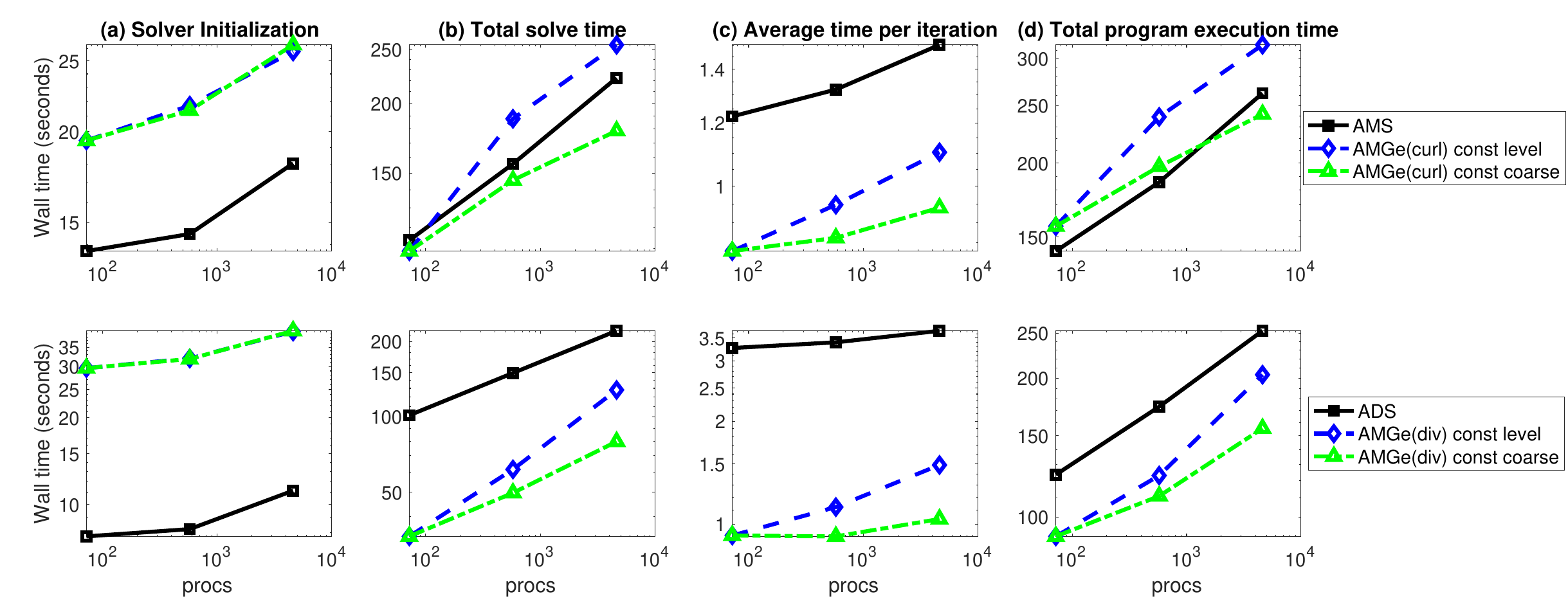}
\caption[]{Solver weak scaling with lowest order finite elements, where elems / procs = 103,680: $H(\ucurl)$-problem (top row) and $H(\div)$-problem (bottom row).}\label{fig:lo}
\end{figure}

Problem information and solvers iterations are shown in \Cref{tbl:lo} for both $H(\ucurl)$ and $H(\div)$. Timing plots, using wall-clock times as reported by the code of the miniapp, are shown in \Cref{fig:lo}, including wall-clock times for the entire program executions. In the legend of \Cref{fig:lo}, AMGe($\ucurl$) and AMGe($\div$) denote the AMGe solvers for the $H(\ucurl)$ and $H(\div)$ problems, respectively. The cases of constant number of levels  and a constant size of the coarsest problem are indicated with \emph{const levels} and \emph{const coarse}. Finally, AMS and ADS represent the solvers from the HYPRE library.

\begin{table}
\centering
\footnotesize
\begin{tabular}{ | c | r | r | c | r | r | r | r | r | r | }
\hline
\multirow{2}{*}{refs} & \multirow{2}{*}{procs} & \multicolumn{1}{|c|}{\multirow{2}{*}{elems}} & \multirow{2}{*}{$\ell$} & \multicolumn{3}{|c|}{$H(\ucurl)$} & \multicolumn{3}{|c|}{$H(\div)$} \\
\cline{5-10}
& & & & \multicolumn{1}{|c|}{dofs} & $\mathrm{it}_\mathrm{e}$ & $\mathrm{it}_\mathrm{a}$ & \multicolumn{1}{|c|}{dofs} & $\mathrm{it}_\mathrm{e}$ & $\mathrm{it}_\mathrm{a}$ \\
\hline
2 & 720 & 7,464,960 & 3 & 180,667,656 & 191 & 166 & 179,912,448 & 60 & 57 \\
\hline
3 & 5,760 & 59,719,680 & 3 & 1,439,303,184 & 269 & 223 & 1,436,285,952 & 100 & 78 \\
\hline
3 & 5,760 & 59,719,680 & 4 & 1,439,303,184 & 285 & 223 & 1,436,285,952 & 84 & 78 \\
\hline
\end{tabular}
\caption[]{
Weak scaling results with next to the lowest order finite elements for both $H(\ucurl)$ and $H(\div)$ problems, as provided by \eqref{eq:bfs} and \Cref{fig:coeff}. Here, $\ell$ denotes the number of levels in the AMGe hierarchy, $\mathrm{it}_\mathrm{e}$ -- the number of PCG iteration using the proposed AMGe preconditioner, and $\mathrm{it}_\mathrm{a}$ -- the number of PCG iteration using the auxiliary space (AMS or ADS) preconditioners from the HYPRE library. In all cases, elems / procs = 10,368, while the GC and OC of the AMGe hierarchy are approximately constant and equal to 1.09 for all refinement levels.
\label{tbl:ho}
}
\end{table}

\begin{figure}
\centering
\includegraphics[width=\textwidth]{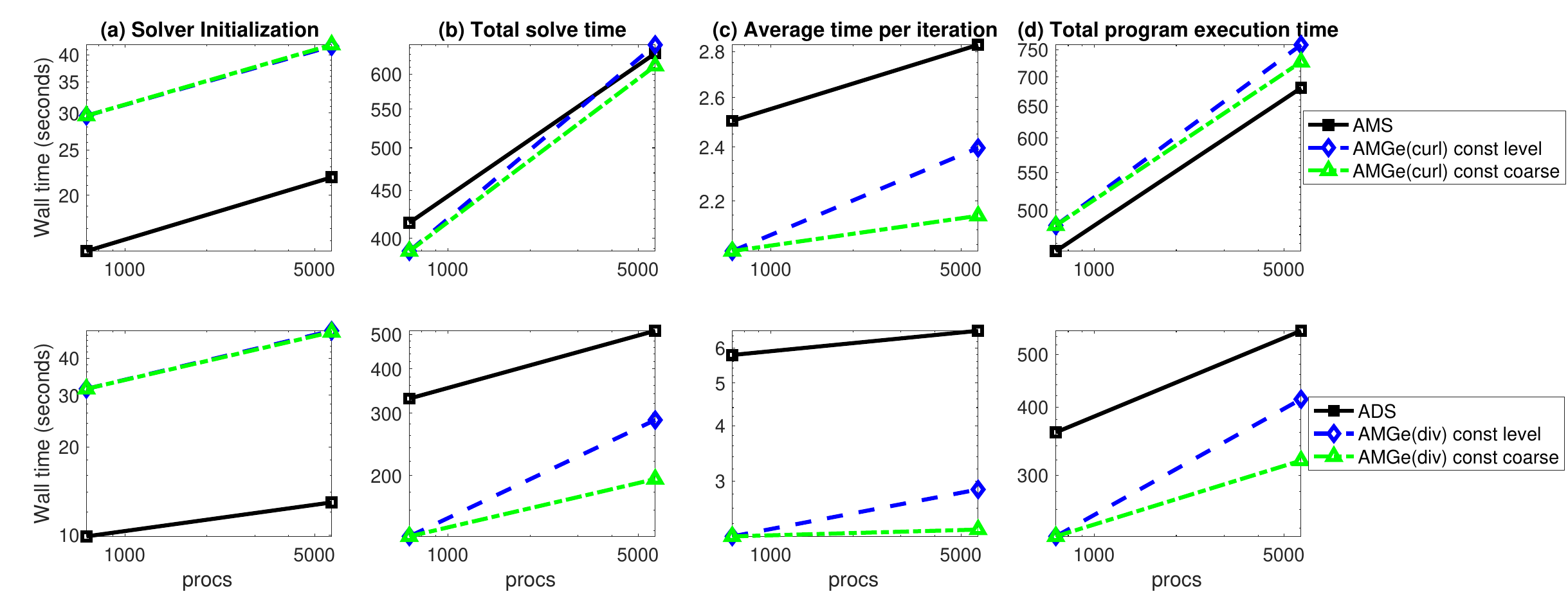}
\caption[]{Solver weak scaling with next to the lowest order finite elements, where elems / procs = 10,368: $H(\ucurl)$-problem (top row) and $H(\div)$-problem (bottom row).}\label{fig:ho}
\end{figure}

The AMGe approach delivers good and competitive performance in comparison with the state of the art represented by AMS and ADS. Moreover, due to the higher set-up cost but shorter solve times, ParELAG has a competitive advantage compared the AMS and ADS when the same hierarchy can be reused in solving multiple linear systems with different right-hand sides. Observe that \Cref{fig:lo} indicates that the more standard case of utilizing AMGe, of increasing the number of AMGe levels, demonstrates better scalability. To further study and exploit this scalability potential in practice for extremely large problems in the setting of extreme parallelism, parallel redistribution and load balancing are needed on coarse levels obtained via AMGe to allow sufficient coarsening when large number of processors are utilized. This is a subject of an ongoing work.

\subsubsection{Results for next to the lowest order elements}

For completeness, results using next to the lowest order finite elements are presented in \Cref{tbl:ho,fig:ho}, following the paradigm of the previous subsection. Again, the AMGe methodology performs well and is comparable to the state of the art.

\begin{remark}
Note that the finest MFEM-generated level and the coarse ParELAG-generated levels in the miniapp by default interpret \emph{next to the lowest order} slightly differently, even if similar, especially when employing hexahedral elements. Consider for example the $L^2$-conforming finite elements spaces of piecewise polynomial functions. Being informed about the geometry of the elements and their tensor-product structure, MFEM produces bilinear elements with 8 dofs and basis functions per element. In contrast, ParELAG operates in a generic geometry-agnostic way and by default the finite element order determines the order of the targets. Thus, by default it produces targets and spaces that provide piecewise linear interpolation independently of the shape of the element, resulting in 4 dofs and basis functions per element. Generally, this behavior can be easily altered by changing the way targets are selected, but the concentration here is on the default behavior.
\end{remark}

\section{Conclusions and future work}
\label{sec:conclusions}

In this paper, we have introduced an AMGe approach for $H(\ucurl)$ and $H(\div)$ formulations. It involves the construction of coarse de Rham sequences on agglomerated meshes, the use of hybrid (Hiptmair) smoothers, and state-of-the-art auxiliary space multigrid solvers, HYPRE AMS and ADS, for the coarsest level. The methods are described in detail using the exterior calculus framework, which allows for a unified (independent of the number of space dimensions) presentation of the local problems that need be solved to construct the hierarchy of de Rham sequences. A key characteristic of the AMGe technique implemented in ParELAG is that the de Rham sequences at each level of the hierarchy possess the same properties and structures that are defined at the finite element levels, including boundary attributes, exterior derivative operators, co-chain projectors, and interpolation operators between spaces. The paper also provide an overview of the ParELAG implementation of the above methods. The numerical results presented here demonstrate the good performance and weak scaling properties of ParELAG, comparable to the state-of-the-art $H(\ucurl)$ and $H(\div)$ solvers in HYPRE. A potential limitation of the current implementation is that ParELAG does not, yet, admit agglomerated elements that are shared or redistributed between processors. The development of such functionality is a currently ongoing work, which would also allow us to exploit the improved scaling of deep AMGe cycles. Future work also includes the implementation of AMGe variants of the hybridization and static condensation techniques in \cite{2019HybridizationHdiv,2020PreconditionerIP,2020PreconditionerMortar} for solving the coarsest $H(\div)$ problems.

\section*{Acknowledgements}
The authors would like to acknowledge all the other developers and contributors to the ParELAG library, including Andrew Barker, Thomas Benson, Ilya Lashuk, Chak Shing Lee, and Sara Osborn (contributors listed in alphabetical order).

\bibliographystyle{plainurl}
\bibliography{references.bib}

\begin{thebibliography}{10}

\bibitem{hypre}
{HYPRE: Scalable Linear Solvers and Multigrid Methods}.
\newblock
  \url{http://computing.llnl.gov/projects/hypre-scalable-linear-solvers-multigrid-methods}.

\bibitem{metis}
{METIS: Graph Partitioning and Fill-reducing Matrix Ordering}.
\newblock \url{http://glaros.dtc.umn.edu/gkhome/views/metis}.

\bibitem{mfem}
{MFEM: Modular Finite Element Methods Library}.
\newblock \url{http://mfem.org}.
\newblock \href {https://doi.org/10.11578/dc.20171025.1248}
  {\path{doi:10.11578/dc.20171025.1248}}.

\bibitem{miniapps}
{ParELAG mini applications in MFEM}.
\newblock \url{http://github.com/mfem/mfem/tree/master/miniapps/parelag}.

\bibitem{parelag}
{ParELAG: Parallel Element Agglomeration Algebraic Multigrid Upscaling and
  Solvers}.
\newblock \url{http://github.com/LLNL/parelag}.

\bibitem{strumpack}
{STRUMPACK: STRUctured Matrix PACKage}.
\newblock \url{http://portal.nersc.gov/project/sparse/strumpack}.

\bibitem{2013ConstrFOSLS}
J~H Adler and P~S Vassilevski.
\newblock {Improving Conservation for First-Order System Least-Squares
  Finite-Element Methods}.
\newblock In Oleg~P Iliev, Svetozar~D Margenov, Peter~D Minev, Panayot~S
  Vassilevski, and Ludmil~T Zikatanov, editors, {\em Numer. Solut. Partial
  Differ. Equations Theory, Algorithms, Their Appl.}, pages 1--19, 2013.
\newblock \href {https://doi.org/10.1007/978-1-4614-7172-1_1}
  {\path{doi:10.1007/978-1-4614-7172-1_1}}.

\bibitem{AdlerVassilevski2014}
JH~Adler and Panayot~S Vassilevski.
\newblock Error analysis for constrained first-order system least-squares
  finite-element methods.
\newblock {\em {SIAM Journal on Scientific Computing}}, 36(3):A1071--A1088,
  2014.

\bibitem{2012VecLaplace}
Douglas~N Arnold, Richard~S Falk, and Jay Gopalakrishnan.
\newblock {Mixed Finite Element Approximation of the Vector Laplace with
  Dirichlet Boundary Conditions}.
\newblock {\em Math. Model. Methods Appl. Sci.}, 22(09):1250024, 2012.
\newblock \href {https://doi.org/10.1142/S0218202512500248}
  {\path{doi:10.1142/S0218202512500248}}.

\bibitem{1997HdivPecond}
Douglas~N. Arnold, Richard~S. Falk, and Ragnar Winther.
\newblock {Preconditioning in $\mathbf{H}(\operatorname{div})$ and
  applications}.
\newblock {\em Math. Comput.}, 66(219):957--985, 1997.
\newblock \href {https://doi.org/10.1090/S0025-5718-97-00826-0}
  {\path{doi:10.1090/S0025-5718-97-00826-0}}.

\bibitem{2000HdivHcurlMultigrid}
Douglas~N Arnold, Richard~S Falk, and Ragnar Winther.
\newblock {Multigrid in $H(\operatorname{div})$ and $H(\operatorname{curl})$}.
\newblock {\em Numer. Math.}, 85(2):197--217, 2000.
\newblock \href {https://doi.org/10.1007/PL00005386}
  {\path{doi:10.1007/PL00005386}}.

\bibitem{2010DeRhamReview}
Douglas~N. Arnold, Richard~S. Falk, and Ragnar Winther.
\newblock {Finite element exterior calculus: from Hodge theory to numerical
  stability}.
\newblock {\em Bull. Am. Math. Soc.}, 47(2):281--354, 2010.
\newblock \href {https://doi.org/10.1090/S0273-0979-10-01278-4}
  {\path{doi:10.1090/S0273-0979-10-01278-4}}.

\bibitem{2011Smoothers}
A~Baker, R~Falgout, T~Kolev, and U~Yang.
\newblock {Multigrid Smoothers for Ultraparallel Computing}.
\newblock {\em SIAM J. Sci. Comput.}, 33(5):2864--2887, 2011.
\newblock \href {https://doi.org/10.1137/100798806}
  {\path{doi:10.1137/100798806}}.

\bibitem{2008AMGkLaplacian}
Nathan Bell and Luke~N Olson.
\newblock {Algebraic multigrid for $k$-form Laplacians}.
\newblock {\em Numer. Linear Algebr. with Appl.}, 15(2‐3):165--185, 2008.
\newblock \href {https://doi.org/10.1002/nla.577} {\path{doi:10.1002/nla.577}}.

\bibitem{2003HcurlAMG}
Pavel~B Bochev, Christopher~J Garasi, Jonathan~J Hu, Allen~C Robinson, and
  Raymond~S Tuminaro.
\newblock {An Improved Algebraic Multigrid Method for Solving Maxwell's
  Equations}.
\newblock {\em SIAM J. Sci. Comput.}, 25(2):623--642, 2003.
\newblock \href {https://doi.org/10.1137/S1064827502407706}
  {\path{doi:10.1137/S1064827502407706}}.

\bibitem{2003MagneticDiffusion}
Pavel~B Bochev, Jonathan~J Hu, Allen~C Robinson, and Raymond~S Tuminaro.
\newblock {Towards robust 3D Z-pinch simulations: Discretization and fast
  solvers for magnetic diffusion in heterogeneous conductors}.
\newblock {\em Electron. Trans. Numer. Anal.}, 15:186--210, 2003.

\bibitem{2008AuxiliarySpace}
Pavel~B Bochev, Jonathan~J Hu, Christopher~M Siefert, and Raymond~S Tuminaro.
\newblock {An Algebraic Multigrid Approach Based on a Compatible Gauge
  Reformulation of Maxwell's Equations}.
\newblock {\em SIAM J. Sci. Comput.}, 31(1):557--583, 2008.
\newblock \href {https://doi.org/10.1137/070685932}
  {\path{doi:10.1137/070685932}}.

\bibitem{BoffiMFE}
Daniele Boffi, Franco Brezzi, and Michel Fortin.
\newblock {\em {Mixed Finite Element Methods and Applications}}, volume~44 of
  {\em Springer Series in Computational Mathematics}.
\newblock Springer, Berlin, Heidelberg, 2013.

\bibitem{2012SAAMGe}
Marian Brezina and Panayot~S Vassilevski.
\newblock {Smoothed Aggregation Spectral Element Agglomeration AMG:
  SA-$\rho$AMGe}.
\newblock In Ivan Lirkov, Svetozar Margenov, and Jerzy Wa{\'{s}}niewski,
  editors, {\em Large-Scale Sci. Comput.}, pages 3--15, Berlin, Heidelberg,
  2012. Springer.

\bibitem{Brunner2002}
Thomas~A Brunner.
\newblock {Forms of Approximate Radiation Transport}.
\newblock Technical report, SAND2002-1778, Sandia National Laboratories, 2002.
\newblock \href {https://doi.org/10.2172/800993} {\path{doi:10.2172/800993}}.

\bibitem{1994FOSLS1}
Z~Cai, R~Lazarov, T~A Manteuffel, and S~F McCormick.
\newblock {First-Order System Least Squares for Second-Order Partial
  Differential Equations: Part I}.
\newblock {\em SIAM J. Numer. Anal.}, 31(6):1785--1799, 1994.
\newblock \href {https://doi.org/10.1137/0731091} {\path{doi:10.1137/0731091}}.

\bibitem{2010MixedStokes}
Zhiqiang Cai, Charles Tong, Panayot~S Vassilevski, and Chunbo Wang.
\newblock {Mixed finite element methods for incompressible flow: Stationary
  Stokes equations}.
\newblock {\em Numer. Methods Partial Differ. Equ.}, 26(4):957--978, 2010.
\newblock \href {https://doi.org/10.1002/num.20467}
  {\path{doi:10.1002/num.20467}}.

\bibitem{2010MixedNS}
Zhiqiang Cai, Chunbo Wang, and Shun Zhang.
\newblock {Mixed Finite Element Methods for Incompressible Flow: Stationary
  Navier-Stokes Equations}.
\newblock {\em SIAM J. Numer. Anal.}, 48(1):79--94, 2010.
\newblock \href {https://doi.org/10.1137/080718413}
  {\path{doi:10.1137/080718413}}.

\bibitem{2003AMGe}
T~Chartier, R~Falgout, V~Henson, J~Jones, T~Manteuffel, S~McCormick, J~Ruge,
  and P~Vassilevski.
\newblock {Spectral AMGe ($\rho$AMGe)}.
\newblock {\em SIAM J. Sci. Comput.}, 25(1):1--26, 2003.
\newblock \href {https://doi.org/10.1137/S106482750139892X}
  {\path{doi:10.1137/S106482750139892X}}.

\bibitem{2019HybridizationHdiv}
V~Dobrev, T~Kolev, C~S Lee, V~Tomov, and P~S Vassilevski.
\newblock {Algebraic Hybridization and Static Condensation with Application to
  Scalable $H(\mathrm{div})$ Preconditioning}.
\newblock {\em SIAM J. Sci. Comput.}, 41(3):B425--B447, 2019.
\newblock \href {https://doi.org/10.1137/17M1132562}
  {\path{doi:10.1137/17M1132562}}.

\bibitem{2021PosteriorMultilevel}
Hillary~R Fairbanks, Sarah Osborn, and Panayot~S Vassilevski.
\newblock {Estimating posterior quantity of interest expectations in a
  multilevel scalable framework}.
\newblock {\em Numer. Linear Algebr. with Appl.}, 28(3):e2352, 2021.
\newblock \href {https://doi.org/10.1002/nla.2352}
  {\path{doi:10.1002/nla.2352}}.

\bibitem{2021HierarchicalMLMCMC}
Hillary~R Fairbanks, Umberto Villa, and Panayot~S Vassilevski.
\newblock {Multilevel Hierarchical Decomposition of Finite Element White Noise
  with Application to Multilevel Markov Chain Monte Carlo}.
\newblock {\em SIAM J. Sci. Comput.}, pages S293--S316, 2021.
\newblock \href {https://doi.org/10.1137/20M1349606}
  {\path{doi:10.1137/20M1349606}}.

\bibitem{2004GenAMG}
Robert~D Falgout and Panayot~S Vassilevski.
\newblock {On Generalizing the Algebraic Multigrid Framework}.
\newblock {\em SIAM J. Numer. Anal.}, 42(4):1669--1693, 2004.
\newblock \href {https://doi.org/10.1137/S0036142903429742}
  {\path{doi:10.1137/S0036142903429742}}.

\bibitem{gentile2001implicit}
NA~Gentile.
\newblock Implicit {M}onte {C}arlo diffusion—an acceleration method for
  {M}onte {C}arlo time-dependent radiative transfer simulations.
\newblock {\em Journal of Computational Physics}, 172(2):543--571, 2001.

\bibitem{golub1999inexact}
Gene~H Golub and Qiang Ye.
\newblock Inexact preconditioned conjugate gradient method with inner-outer
  iteration.
\newblock {\em SIAM Journal on Scientific Computing}, 21(4):1305--1320, 1999.

\bibitem{2018AuxiliaryDeRham}
J~Gopalakrishnan, M~Neum{\"{u}}ller, and P~S Vassilevski.
\newblock {The Auxiliary Space Preconditioner for the de Rham Complex}.
\newblock {\em SIAM J. Numer. Anal.}, 56(6):3196--3218, 2018.
\newblock \href {https://doi.org/10.1137/17M1153376}
  {\path{doi:10.1137/17M1153376}}.

\bibitem{graziani2000crooked}
F~Graziani and J~LeBlanc.
\newblock The crooked pipe test problem.
\newblock {\em Lawrence Livermore National Laboratory Report UCRL-MI-143393},
  2000.

\bibitem{1997Hiptmair}
R~Hiptmair.
\newblock {Multigrid method for $H(\operatorname{div})$ in three dimensions}.
\newblock {\em Electron. Trans. Numer. Anal.}, 6:133--152, 1997.

\bibitem{1998Hiptmair}
R~Hiptmair.
\newblock {Multigrid Method for Maxwell's Equations}.
\newblock {\em SIAM J. Numer. Anal.}, 36(1):204--225, 1998.
\newblock \href {https://doi.org/10.1137/S0036142997326203}
  {\path{doi:10.1137/S0036142997326203}}.

\bibitem{2002FEelectromagDeRham}
R~Hiptmair.
\newblock {Finite elements in computational electromagnetism}.
\newblock {\em Acta Numer.}, 11:237--339, 2002.
\newblock \href {https://doi.org/10.1017/S0962492902000041}
  {\path{doi:10.1017/S0962492902000041}}.

\bibitem{2006HcurlAux}
R~Hiptmair, G~Widmer, and J~Zou.
\newblock {Auxiliary space preconditioning in $H_0(\mathrm{curl}; \Omega)$}.
\newblock {\em Numer. Math.}, 103(3):435--459, 2006.
\newblock \href {https://doi.org/10.1007/s00211-006-0683-0}
  {\path{doi:10.1007/s00211-006-0683-0}}.

\bibitem{2007AuxiliarySpace}
R~Hiptmair and J~Xu.
\newblock {Nodal Auxiliary Space Preconditioning in $\mathbf{H}(\mathbf{curl})$
  and $\mathbf{H}(\operatorname{div})$ Spaces}.
\newblock {\em SIAM J. Numer. Anal.}, 45(6):2483--2509, 2007.
\newblock \href {https://doi.org/10.1137/060660588}
  {\path{doi:10.1137/060660588}}.

\bibitem{2000SchwarzHcurlHdiv}
Ralf Hiptmair and Andrea Toselli.
\newblock {Overlapping and Multilevel Schwarz Methods for Vector Valued
  Elliptic Problems in Three Dimensions}.
\newblock In Petter Bj{\o}rstad and Mitchell Luskin, editors, {\em Parallel
  Solut. Partial Differ. Equations}, pages 181--208, 2000.
\newblock \href {https://doi.org/10.1007/978-1-4612-1176-1_8}
  {\path{doi:10.1007/978-1-4612-1176-1_8}}.

\bibitem{2006HcurlAMG}
J~Jones and B~Lee.
\newblock {A Multigrid Method for Variable Coefficient Maxwell's Equations}.
\newblock {\em SIAM J. Sci. Comput.}, 27(5):1689--1708, 2006.
\newblock \href {https://doi.org/10.1137/040608283}
  {\path{doi:10.1137/040608283}}.

\bibitem{2001AMGe}
Jim~E Jones and Panayot~S Vassilevski.
\newblock {AMGe Based on Element Agglomeration}.
\newblock {\em SIAM J. Sci. Comput.}, 23(1):109--133, 2001.
\newblock \href {https://doi.org/10.1137/S1064827599361047}
  {\path{doi:10.1137/S1064827599361047}}.

\bibitem{2016AMGeUpscaling}
D~Z Kalchev, C~S Lee, U~Villa, Y~Efendiev, and P~S Vassilevski.
\newblock {Upscaling of Mixed Finite Element Discretization Problems by the
  Spectral AMGe Method}.
\newblock {\em SIAM J. Sci. Comput.}, 38(5):A2912--A2933, 2016.
\newblock \href {https://doi.org/10.1137/15M1036683}
  {\path{doi:10.1137/15M1036683}}.

\bibitem{2020PreconditionerMortar}
Delyan~Z Kalchev and Panayot Vassilevski.
\newblock {A Condensed Constrained Nonconforming Mortar-Based Approach for
  Preconditioning Finite Element Discretization Problems}.
\newblock {\em SIAM J. Sci. Comput.}, 42(5):A3136--A3156, 2020.
\newblock \href {https://doi.org/10.1137/19M1305690}
  {\path{doi:10.1137/19M1305690}}.

\bibitem{2020PreconditionerIP}
Delyan~Z Kalchev and Panayot~S Vassilevski.
\newblock {Auxiliary Space Preconditioning of Finite Element Equations Using a
  Nonconforming Interior Penalty Reformulation and Static Condensation}.
\newblock {\em SIAM J. Sci. Comput.}, 42(3):A1741--A1764, 2020.
\newblock \href {https://doi.org/10.1137/19M1286815}
  {\path{doi:10.1137/19M1286815}}.

\bibitem{2008HcurlAux}
Tzanio~V Kolev, Joseph~E Pasciak, and Panayot~S Vassilevski.
\newblock {$\mathbf{H}(\mathrm{curl})$ auxiliary mesh preconditioning}.
\newblock {\em Numer. Linear Algebr. with Appl.}, 15(5):455--471, 2008.
\newblock \href {https://doi.org/10.1002/nla.534} {\path{doi:10.1002/nla.534}}.

\bibitem{2009AMS}
Tzanio~V Kolev and Panayot~S Vassilevski.
\newblock {Parallel Auxiliary Space AMG for $H(\operatorname{curl})$ Problems}.
\newblock {\em J. Comput. Math.}, 27(5):604--623, 2009.
\newblock \href {https://doi.org/10.4208/jcm.2009.27.5.013}
  {\path{doi:10.4208/jcm.2009.27.5.013}}.

\bibitem{2012ADS}
Tzanio~V Kolev and Panayot~S Vassilevski.
\newblock {Parallel Auxiliary Space AMG Solver for $H(\operatorname{div})$
  Problems}.
\newblock {\em SIAM J. Sci. Comput.}, 34(6):A3079--A3098, 2012.
\newblock \href {https://doi.org/10.1137/110859361}
  {\path{doi:10.1137/110859361}}.

\bibitem{2018AMGeFAS}
Max {la Cour Christensen}, Panayot~S Vassilevski, and Umberto Villa.
\newblock {Nonlinear multigrid solvers exploiting AMGe coarse spaces with
  approximation properties}.
\newblock {\em J. Comput. Appl. Math.}, 340:691--708, 2018.
\newblock \href {https://doi.org/10.1016/j.cam.2017.10.029}
  {\path{doi:10.1016/j.cam.2017.10.029}}.

\bibitem{2017ParELAGReservoir}
Max {la Cour Christensen}, Umberto Villa, Allan~P Engsig-Karup, and Panayot~S
  Vassilevski.
\newblock {Numerical Multilevel Upscaling for Incompressible Flow in Reservoir
  Simulation: An Element-Based Algebraic Multigrid (AMGe) Approach}.
\newblock {\em SIAM J. Sci. Comput.}, 39(1):B102--B137, 2017.
\newblock \href {https://doi.org/10.1137/140988991}
  {\path{doi:10.1137/140988991}}.

\bibitem{2012AMGeRT}
I~V Lashuk and P~S Vassilevski.
\newblock {Element agglomeration coarse Raviart-Thomas spaces with improved
  approximation properties}.
\newblock {\em Numer. Linear Algebr. with Appl.}, 19(2):414--426, 2012.
\newblock \href {https://doi.org/10.1002/nla.1819}
  {\path{doi:10.1002/nla.1819}}.

\bibitem{2008AMGe}
Ilya Lashuk and Panayot~S Vassilevski.
\newblock {On some versions of the element agglomeration AMGe method}.
\newblock {\em Numer. Linear Algebr. with Appl.}, 15(7):595--620, 2008.
\newblock \href {https://doi.org/10.1002/nla.585} {\path{doi:10.1002/nla.585}}.

\bibitem{2014CoarseDeRham}
Ilya~V Lashuk and Panayot~S Vassilevski.
\newblock {The Construction of the Coarse de Rham Complexes with Improved
  Approximation Properties}.
\newblock {\em Comput. Methods Appl. Math.}, 14(2):257--303, 2014.
\newblock \href {https://doi.org/10.1515/cmam-2014-0004}
  {\path{doi:10.1515/cmam-2014-0004}}.

\bibitem{superludist}
Xiaoye~S Li and James~W Demmel.
\newblock {SuperLU\_DIST: A Scalable Distributed-Memory Sparse Direct Solver
  for Unsymmetric Linear Systems}.
\newblock {\em ACM Trans. Math. Softw.}, 29(2):110--140, jun 2003.
\newblock \href {https://doi.org/10.1145/779359.779361}
  {\path{doi:10.1145/779359.779361}}.

\bibitem{1997NS}
Ping Lin.
\newblock {A Sequential Regularization Method for Time-Dependent Incompressible
  Navier-Stokes Equations}.
\newblock {\em SIAM J. Numer. Anal.}, 34(3):1051--1071, 1997.
\newblock \href {https://doi.org/10.1137/S0036142994270521}
  {\path{doi:10.1137/S0036142994270521}}.

\bibitem{MonkFEM}
Peter Monk.
\newblock {\em {Finite Element Methods for Maxwell's Equations}}.
\newblock Numerical Mathematics and Scientific Computation. Clarendon Press,
  Oxford, 2003.

\bibitem{2018HdivBDDC}
Duk-Soon Oh, Olof~B Widlund, Stefano Zampini, and Clark~R Dohrmann.
\newblock {BDDC Algorithms with deluxe scaling and adaptive selection of primal
  constraints for Raviart-Thomas vector fields}.
\newblock {\em Math. Comput.}, 87(310):659--692, 2018.
\newblock \href {https://doi.org/10.1090/mcom/3254}
  {\path{doi:10.1090/mcom/3254}}.

\bibitem{2017PDESampler}
Sarah Osborn, Panayot~S Vassilevski, and Umberto Villa.
\newblock {A Multilevel, Hierarchical Sampling Technique for Spatially
  Correlated Random Fields}.
\newblock {\em SIAM J. Sci. Comput.}, 39(5):S543--S562, 2017.
\newblock \href {https://doi.org/10.1137/16M1082688}
  {\path{doi:10.1137/16M1082688}}.

\bibitem{2018PDESampler}
Sarah Osborn, Patrick Zulian, Thomas Benson, Umberto Villa, Rolf Krause, and
  Panayot~S Vassilevski.
\newblock {Scalable hierarchical PDE sampler for generating spatially
  correlated random fields using nonmatching meshes}.
\newblock {\em Numer. Linear Algebr. with Appl.}, 25(3):e2146, 2018.
\newblock \href {https://doi.org/10.1002/nla.2146}
  {\path{doi:10.1002/nla.2146}}.

\bibitem{2002HcurlSchwartz}
J~E Pasciak and J~Zhao.
\newblock {Overlapping Schwarz methods in $\mathbf{H}(\mathbf{curl})$ on
  polyhedral domains}.
\newblock {\em J. Numer. Math.}, 10(3):221--234, 2002.
\newblock \href {https://doi.org/10.1515/JNMA.2002.221}
  {\path{doi:10.1515/JNMA.2002.221}}.

\bibitem{2008deRhamAMGe}
Joseph~E Pasciak and Panayot~S Vassilevski.
\newblock {Exact de Rham Sequences of Spaces Defined on Macro-Elements in Two
  and Three Spatial Dimensions}.
\newblock {\em SIAM J. Sci. Comput.}, 30(5):2427--2446, 2008.
\newblock \href {https://doi.org/10.1137/070698178}
  {\path{doi:10.1137/070698178}}.

\bibitem{1996LSnonSA}
A~I Pehlivanov, G~F Carey, and P~S Vassilevski.
\newblock {Least-squares mixed finite element methods for non-selfadjoint
  elliptic problems: I. Error estimates}.
\newblock {\em Numer. Math.}, 72(4):501--522, 1996.
\newblock \href {https://doi.org/10.1007/s002110050179}
  {\path{doi:10.1007/s002110050179}}.

\bibitem{2002HcurlAMG}
S~Reitzinger and J~Sch{\"{o}}berl.
\newblock {An algebraic multigrid method for finite element discretizations
  with edge elements}.
\newblock {\em Numer. Linear Algebr. with Appl.}, 9(3):223--238, 2002.
\newblock \href {https://doi.org/10.1002/nla.271} {\path{doi:10.1002/nla.271}}.

\bibitem{2007LEforMHD}
R~N Rieben, D~A White, B~K Wallin, and J~M Solberg.
\newblock {An arbitrary Lagrangian-Eulerian discretization of MHD on 3D
  unstructured grids}.
\newblock {\em J. Comput. Phys.}, 226(1):534--570, 2007.
\newblock \href {https://doi.org/10.1016/j.jcp.2007.04.031}
  {\path{doi:10.1016/j.jcp.2007.04.031}}.

\bibitem{vass_sparse_matrix_topology}
Panayot~S. Vassilevski.
\newblock Sparse matrix element topology with application to amg and
  preconditioning.
\newblock {\em Numer. Lin. Alg. Appl.}, 9:429--444, 2002.

\bibitem{VassilevskiMG}
Panayot~S Vassilevski.
\newblock {\em {Multilevel Block Factorization Preconditioners: Matrix-based
  Analysis and Algorithms for Solving Finite Element Equations}}.
\newblock Springer, New York, 2008.
\newblock \href {https://doi.org/10.1007/978-0-387-71564-3}
  {\path{doi:10.1007/978-0-387-71564-3}}.

\bibitem{2011UpscalingAMG}
Panayot~S Vassilevski.
\newblock {Coarse Spaces by Algebraic Multigrid: Multigrid Convergence and
  Upscaling Error Estimates}.
\newblock {\em Adv. Adapt. Data Anal.}, 03(01n02):229--249, 2011.
\newblock \href {https://doi.org/10.1142/S1793536911000830}
  {\path{doi:10.1142/S1793536911000830}}.

\bibitem{2013Brinkman}
Panayot~S Vassilevski and Umberto Villa.
\newblock {A Block-Diagonal Algebraic Multigrid Preconditioner for the Brinkman
  Problem}.
\newblock {\em SIAM J. Sci. Comput.}, 35(5):S3--S17, 2013.
\newblock \href {https://doi.org/10.1137/120882846}
  {\path{doi:10.1137/120882846}}.

\bibitem{1992MixedMultigrid}
Panayot~S. Vassilevski and Junping Wang.
\newblock {Multilevel iterative methods for mixed finite element
  discretizations of elliptic problems}.
\newblock {\em Numer. Math.}, 63(1):503--520, 1992.
\newblock \href {https://doi.org/10.1007/BF01385872}
  {\path{doi:10.1007/BF01385872}}.

\bibitem{2009OptimalMG}
Jinchao Xu, Long Chen, and Ricardo~H Nochetto.
\newblock {Optimal multilevel methods for $H(\mathrm{grad})$,
  $H(\mathrm{curl})$, and $H(\mathrm{div})$ systems on graded and unstructured
  grids}.
\newblock In Ronald DeVore and Angela Kunoth, editors, {\em Multiscale,
  Nonlinear Adapt. Approx.}, pages 599--659, Berlin, Heidelberg, 2009.
  Springer.
\newblock \href {https://doi.org/10.1007/978-3-642-03413-8_14}
  {\path{doi:10.1007/978-3-642-03413-8_14}}.

\bibitem{2016PETScBDDC}
Stefano Zampini.
\newblock {PCBDDC: A Class of Robust Dual-Primal Methods in PETSc}.
\newblock {\em SIAM J. Sci. Comput.}, 38(5):S282--S306, 2016.
\newblock \href {https://doi.org/10.1137/15M1025785}
  {\path{doi:10.1137/15M1025785}}.

\bibitem{2016RobustnessBDDC}
Stefano Zampini and David~E Keyes.
\newblock {On the Robustness and Prospects of Adaptive BDDC Methods for Finite
  Element Discretizations of Elliptic PDEs with High-Contrast Coefficients}.
\newblock In {\em Proc. Platf. Adv. Sci. Comput. Conf.}, New York, 2016.
  Association for Computing Machinery.
\newblock \href {https://doi.org/10.1145/2929908.2929919}
  {\path{doi:10.1145/2929908.2929919}}.

\end{thebibliography}

\end{document}